\documentclass[10pt,a4paper,twoside]{article}
\usepackage{epsfig,color,graphics,graphicx}
\usepackage{amssymb,amsbsy,amsmath,amsfonts,amssymb,amscd}
\usepackage[english]{babel}
\newtheorem{theorem}{Theorem}

\newtheorem{corollary}[theorem]{Corollary}

\newtheorem{lemma}[theorem]{Lemma}

\newtheorem{proposition}[theorem]{Proposition}
\newtheorem{remark}{Remark}

\newenvironment{proof}[1][Proof]{\noindent\textbf{#1.} }{\ \rule{0.5em}{0.5em}}

\newcommand{\m}{{\mathbf m}}
\newcommand{\x}{{\mathbf x}}
\newcommand{\y}{{\mathbf y}}
\newcommand{\bn}{{\boldsymbol \nu}}
\newcommand{\bt}{{\boldsymbol \tau}}
\newcommand{\R}{{\mathbb R}}

\begin{document}

\title{Nonlinear Neumann boundary stabilization\\
of the wave equation using rotated multipliers.}
\author{
{\sc Pierre CORNILLEAU}\thanks{Universit\'e de Lyon, \'Ecole centrale de Lyon,
D.M.I., Institut Camille-Jordan (C.N.R.S. U.M.R. 5208),
36 avenue Guy-de-Collongue, 69134 \'Ecully cedex, France.} ,
{\sc Jean-Pierre LOH\'EAC}$^*$\thanks{Independent Moscow University,
Laboratoire J.-V. Poncelet (C.N.R.S. U.M.I. 2615),
Bol. Vlasyevsky Per. 11, 119002 Moscow, Russia.} ,
{\sc Axel OSSES}\thanks{Universidad de Chile, D.I.M.,
Centro de Modelamiento Matem\'atico, (C.N.R.S. U.M.I. 2807 CNRS-Uchile),
Casilla 170/3--Correo 3, Santiago, Chile.} . 
}
\date{}

\maketitle


\medskip

\begin{abstract}
We study the boundary stabilization of the wave equation by means of a linear
or non-linear Neumann feedback. The rotated multiplier method leads to new
geometrical cases concerning the active part of the boundary where the
feedback is applied. Due to mixed boundary conditions, these cases generate
singularities. Under a simple geometrical condition concerning the orientation
of the boundary, we obtain stabilization results in both cases.
\end{abstract}

\section*{Introduction}
In this paper we are concerned with the stabilization of the wave
equation in a multi-dimensional body $\Omega\subset \R^n $
by using a feedback law applied on some part of its boundary. The problem can
be written as follows
$$ \left\{ 
\begin{matrix}
u''-\Delta u=0\quad
&\text{in } \Omega \times \mathbb{R}^*_+ \, ,\hfill \\
u=0\hfill
&\text{on } \partial \Omega_D \times \mathbb{R}^*_+ \, ,\hfill \\
\partial_\nu u=F \hfill
&\text{on } \partial \Omega_N \times \mathbb{R}^*_+ \, ,\hfill \\
u(0)=u_0 \hfill
&\text{in } \Omega \, ,\hfill \\
u'(0)=u_1 \hfill
&\text{in } \Omega \, ,\hfill
\end{matrix}
\right. $$
where we denote by $u' $, $u'' $, $\Delta u$ and $\partial_\nu u$ the first 
time-derivative
of $u$, the second time-derivative of the scalar function $u$, the standard 
Laplacian of $u$ and the normal outward derivative
of $u$ on $\partial \Omega $, respectively; 
$(\partial \Omega_D $, $\partial \Omega_N )$ is a partition of 
$\partial \Omega $ and $F$ is the feedback function which may depend on the
state $(u,u')$, the position $\x $ and time $t$.\\
Our purpose here is to choose the feedback function $F$ and the active 
part of the boundary, $\partial \Omega_N $, so
that for every initial data, the energy function
$$E(u,t)=\frac{1}{2}\int_\Omega (|u'(t)|^2 +|\nabla u(t)|^2 )\, d\x  \, ,$$
is decreasing with respect to time $t$, and vanishes as $t\longrightarrow
\infty $.\\
Formally, we can write the time derivative of $E$ as follows
$$E'(u,t)=\int_{\partial \Omega_N } F u' \, d\sigma \, ,$$
and a sufficient condition so that $E$ is non-increasing is\, 
$F u'\le 0$ on $\partial \Omega_N $.\\
In the two-dimensional case and 
in the framework of Hilbert Uniqueness Method \cite{L}, it can be shown that 
the energy function is uniformly decreasing as $t$ tends to $\infty $, 
by choosing \, 
$\m : \x \mapsto \x -\x_0$, where $\x_0$ is some given point in $\R^n $ and
$$\partial \Omega_N =\left\{ \x \in \partial \Omega \, / \, 
m (\x ).\nu (\x )>0\, \right\} \, ,\quad
F=-\m .\bn \, u'\, ,$$
where $\bn $ is the normal unit vector pointing outward of $\Omega $.
This method has been performed by many authors, see for instance \cite{KZ}
and references therein. 
Here we extend the above result for rotated multipliers \cite{O1,O2} by
following \cite{BLM}, i.e. we take in account singularities which can appear
when changing boundary conditions along the interface 
$\Gamma =\overline{\partial \Omega_N } \cap \overline{\partial \Omega_D }$.

\section{Notations and main results}
Let $\Omega $ be a bounded open connected set of $\mathbb{R}^{n}(n\geq 2)$ 
such that 
\begin{equation}\label{G}
\partial \Omega \text{ is of class }\mathcal{C}^{2} \text{ in the sense of 
Ne\v{c}as \cite{Ne}}.
\end{equation}
Let $\x_{0}$ be a fixed
point in $\mathbb{R}^{n}$. We denote by $I$ the $n\times n$ identity matrix,
by $A$ a real $n\times n$ skew-symmetric matrix and by $d$ a positive real
number such that
$d^{2}+\Vert A\Vert ^{2}=1$,where $\| \cdot \| $ stands for the usual $2$-norm 
of matrices. 
We now define the following vector function,
\begin{equation*}
\forall \x\in \mathbb{R}^{n}, \quad \m (\x )=(dI+A)(\x -\x_{0} )\, .
\end{equation*}
We consider a partition $(\partial \Omega_N, \partial \Omega_D)$ of 
$\partial \Omega $ such that
\begin{equation}\label{R}
\left \vert 
\begin{array}{l}
\Gamma =\overline{\partial \Omega _{D}}\cap \overline{\partial \Omega _{N}} 
\text{\ is a }\mathcal{C}^{3}\text{-manifold of dimension }n-2,\\
\m .\bn =0 \text{ on } \Gamma,\\
\partial \Omega \cap \omega \ \text{\ is a }\mathcal{C}^{3}\text{-manifold
of dimension }n-1,\\
\mathcal{H}^{n-1}(\partial \Omega_{D})>0,
\end{array}
\right.
\end{equation}
where $\omega$ is a suitable neighborhood of $\Gamma$ and $\mathcal{H}^{n-1}$ 
denotes the usual $(n-1)$-dimensional Hausdorff measure.

\noindent Let $g:$ $\mathbb{R}\rightarrow\mathbb{R}$ be a measurable function 
such that
\begin{equation}\label{F1}
g \text{ is non-decreasing} \quad \text{and}\quad 
\exists K>0\, :\quad |g(s)|\leq K|s|\, 
\text{ a.e..}
\end{equation}
Let us now consider the following wave problem
\begin{equation*}
(S)\, \left\{ 
\begin{array}{l}
u^{\prime \prime }-\Delta u=0 \\ 
u=0 \\ 
\partial _{\nu }u=-\m .\bn \, g(u^{\prime }) \\ 
u(0)=u_{0} \\ 
u^{^{\prime }}(0)=u_{1}
\end{array}
\right. \left. 
\begin{array}{l}
\text{ in }\Omega \times \mathbb{R}_{+}^{\ast }, \\ 
\text{ on }\partial \Omega _{D}\times \mathbb{R}_{+}^{\ast }, \\ 
\text{ on }\partial \Omega _{N}\times \mathbb{R}_{+}^{\ast }, \\ 
\text{ in }\Omega , \\ 
\text{ in }\Omega ,
\end{array}
\right.
\end{equation*}
for some initial data
\begin{equation*}
(u_{0},u_{1})\in \mathrm{H}_{D}^{1}(\Omega )\times \mathrm{L}^{2}(\Omega )
:=\mathrm{H}
\end{equation*}
with $\mathrm{H}_{D}^{1}(\Omega )=\{v\in \mathrm{H}^{1}(\Omega )\, / \, v=0 
\text{ on }\partial \Omega_{D} \}$.

\noindent This problem is well-posed in $\mathrm{H}$.
Indeed, following Komornik \cite{Ko}, we define the non-linear operator 
$\mathcal{A}$ by 
\begin{equation*}
\left.
\begin{array}{l}
\mathcal{A}(u,v)=(-v,-\Delta u), \\ 
D(\mathcal{A})=\{(u,v)\in \mathrm{H}_{D}^{1}(\Omega ) \times 
\mathrm{H}_{D}^{1}(\Omega )\, /\, \Delta u\in \mathrm{L}^{2}(\Omega )
\text{\, and \, }\partial _{\nu }u=-\m .\bn \, g(v)\,
\text{ on }\partial \Omega _{N}\},
\end{array}
\right.
\end{equation*}
so that $(S)$ can be written in the form
\begin{equation*}
\left\{ 
\begin{array}{l}
(u,v)^{\prime }+\mathcal{A}(u,v)=0 \, ,\\ 
(u,v)(0)=(u_{0},u_{1})\, .
\end{array}
\right.
\end{equation*}
It is a classical fact that $\mathcal{A}$ is a maximal-monotone 
operator on $\mathrm{H}$ and that $D(\mathcal{A})$ is dense in $\mathrm{H}$
for the usual norm
(see for instance \cite{B}). Hence, for any initial data $
(u_{0},v_{0})\in $ $D(\mathcal{A})$ there is a unique strong solution $(u,v)$
such that $u\in \mathrm{W}^{1,\infty }(\mathbb{R};\mathrm{H}_{D}^{1}(\Omega ))$
and 
$\Delta u\in \mathrm{L}^{\infty }(\mathbb{R}_{+};\mathrm{L}^{2}(\Omega ))$.
Moreover, for two initial data, the corresponding solutions satisfy
\begin{equation*}
\forall t\geq 0\, ,\quad 
\Vert (u^{1}(t),v^{1}(t))-(u^{2}(t),v^{2}(t))\Vert_{\mathrm{H}}\leq 
C\Vert (u^{1}_{0},v^{1}_{0})-(u^{2}_{0},v^{2}_{0})\Vert_{\mathrm{H}}\, .
\end{equation*}
Using the density of $D(\mathcal{A})$, one can extend the map
\begin{eqnarray*}
      D (\mathcal{A}) & \longrightarrow & \mathrm{H} \\
     (u_o, v_0) & \longmapsto & (u (t), v (t))
\end{eqnarray*}
to a strongly continuous semi-group of contractions $(S(t))_{t\geq 0}$ and 
define for $(u_{0},v_{0})\in \mathrm{H}$ the weak solution 
$(u(t),u'(t))=S(t)(u_{0},u_{1})$ with the regularity 
$u\in \mathcal{C}(\mathbb{R}_{+};\mathrm{H}_{D}^{1}(\Omega ))\cap 
\mathcal{C}^{1}(\mathbb{R}_{+}; \mathrm{L}^{2}(\Omega ))$.
We hence define the energy function of solutions by
$$E(u,0)=\frac{1}{2}\int_{\Omega }(|u_{1}|^{2}+|\nabla u_{0}|^{2})\, d\x , \qquad 
E(u,t)=\frac{1}{2}\int_{\Omega }(|u^{\prime }(t)|^{2}+|\nabla u(t)|^{2})\, d\x 
\quad \text{if }t>0\, .$$
In order to get stabilization results, we need further assumptions
concerning the feedback function $g$
\begin{equation}\label{F2}
\exists p\geq 1 \, ,\ \exists k>0\, ,\quad
|g(s)|\ge k \min \{ |s|,|s|^p \}\, ,\ \hbox{a.e.} \, .
\end{equation}
Concerning the boundary we assume
\begin{equation}\label{S1}
\partial \Omega_{N} \subset
\{\x\in \partial \Omega \, /\, \m (\x ).\bn (\x )\geq 0\}\, ,\quad
\partial \Omega_{D} \subset
\{\x\in \partial \Omega \, /\, \m (\x ).\bn (\x )\leq0\},
\end{equation}
and the additional geometric assumption
\begin{equation}\label{S2}
\m .\bt \le 0 \, \quad \text{on } \Gamma \, ,
\end{equation}
where $\bt (\x )$ is the normal unit vector pointing outward of $\partial
\Omega_N $ at a point $\x \in \Gamma $ when considering $\partial \Omega_N $
as a sub-manifold of $\partial \Omega $.

\begin{remark} It is worth observing that it is not necessary to assume that 
$$\mathcal{H}^{n-1}(\{\x \in \partial \Omega_{N} \, / \, 
\m (\x ).\bn (\x )>0\})>0.$$ 
to get stabilization. In fact, our choice of $\m $ implies such properties 
(see examples in Section 4) whether the energy tends to zero.
\end{remark}

\noindent Since the pioneering work \cite{Ho}, it is now a well-known fact 
that Rellich type relations \cite{Re} are very useful
for the study of control and stabilization of the wave problem.
As we said before, Komornik and Zuazua \cite{KZ} have
shown how these relations can also help us to stabilize the wave problem.
In order to generalize it in higher dimension than 3, the key-problem is
to show the existence of a decomposition of the solution in regular and 
singular parts 
\cite{Gr,KMR} which can be applied
to stabilization problems or control problems. The
first results towards this direction are due to Moussaoui \cite{Mo}, and
Bey-Loh\'{e}ac-Moussaoui \cite{BLM} who also have established a Rellich type
relation in any dimension.\\
In this new case of Neumann feedback deduced from
\cite{O1,O2}, our goal is to generalize those
Rellich relations to get stabilization results about $(S)$
under assumptions \eqref{S1}, \eqref{S2}.\\
As well as in \cite{Ko}, we shall prove here two results of uniform
boundary stabilization.

\subsection*{Exponential boundary stabilization}

We here consider the case when $p=1$ in \eqref{F2}. This is satisfied when
$g$ is linear,
$$\exists \alpha >0 \, : \qquad \forall s\in \mathbb{R} \, ,\quad g(s)=\alpha s\, .$$
In these cases, the energy function is exponentially decreasing.

\begin{theorem}\label{T1}
Assume that geometrical conditions \eqref{R}, \eqref{S1} hold and that the feedback function $g$
satisfies \eqref{F1} and \eqref{F2} with $p=1$.\\
Then under the further geometrical assumption \eqref{S2}, there exist $C>0$
and $T>0$ (independent of $d$)
such that for all initial data in $\mathrm{H}$, the energy of the
solution $u$ of  satisfies
\begin{equation*}
\forall t>\frac{T}{d}\, ,\quad
E(u,t)\le E(u,0)\, \exp \Bigl( 1-\frac{d}{C}t \Bigr) \, .
\end{equation*}
The above constants $C$ and $T$ depend only on the geometry.
\end{theorem}

\subsection*{Rational boundary stabilization}

We here consider the general case and we get rational boundary stabilization.

\begin{theorem}\label{T2}
Assume that geometrical conditions \eqref{R}, \eqref{S1} hold and that the feedback function $g$
satisfies \eqref{F1} and \eqref{F2} with $p>1$.\\
Then under the further geometrical assumption \eqref{S2}, there exist $C>0$
and $T>0$ (independent of $d$)
such that for all initial data in $\mathrm{H}$, the energy of the
solution $u$ of  satisfies
\begin{equation*}
\forall t>\frac{T}{d}\, ,\quad
E(u,t)\le C\, t^{2/(p-1)} \, .
\end{equation*}
where $C$ depends on the initial energy $E(u,0)$.
\end{theorem}

\begin{remark}
Taking advantage of the works of Banasiak-Roach \cite{BR} who generalized
Grisvard's results \cite{Gr} in the piecewise regular case, we will see that 
Theorems \ref{T1} and \ref{T2} remain true in the bi-dimensional case when
assumption \eqref{G} is replaced by following one
\begin{equation}\label{G'}
\begin{matrix}
&\partial \Omega \text{ is a curvilinear polygon of class } {\cal C}^2 \, ,
\hfill \\
& \text{each component of } \partial \Omega \setminus \Gamma \text{ is a }
{\cal C}^2 \hbox{-manifold of dimension } 1\, ,\\
\end{matrix}
\end{equation}
and when condition \eqref{S2}  is replaced by
\begin{equation}\label{S2'}
\begin{matrix}
\forall \x \in \Gamma \, ,\quad  0\le \varpi_\x \le \pi \quad \text{and\quad if }
\varpi_\x =\pi \, ,\, \, \m (\x ).\bt (\x )\le 0 \, .
\end{matrix}
\end{equation}
where $\varpi_\x $ is the angle at the boundary in the point $\x$.
\end{remark}

\noindent These two results are obtained by estimating some integral of the
energy function as well as in \cite{Ko}. This specific estimates are obtained
thanks to an adapted Rellich relation.
\\Hence, this paper is composed of two sections. In the first one we build
convenient Rellich relations and in the second one we use it to prove 
Theorems 1 and 2.

\section{Rellich relations}

\subsection{A regular case}

We can easily build a Rellich relation corresponding to the above vector
field $m $ when considered functions are smooth enough.

\begin{proposition}\label{RR}
Assume that $\Omega $ is an open set of $\mathbb{R}^{n}$ with
boundary of class $\mathcal{C}^{2}$ in the sense of Ne\v{c}as.
If $u$ belongs to $\mathrm{H}^2 (\Omega )$ then
$$2\int_\Omega \Delta u\, (\m .\nabla u)\, d\x  =
d\, (n-2) \int_\Omega |\nabla u|^2 \, d\x  
+\int_{\partial \Omega } \bigl( 2\partial_\nu u\, (\m .\nabla u)
-\m .\bn \, |\nabla u|^2 \bigr) \, d\sigma \, .$$
\end{proposition}

\begin{proof}
Using Green-Riemann identity we get%
\begin{equation*}
2\int_{\Omega }\Delta u\, (\m .\nabla u)\, d\x 
=\int_{\partial \Omega }2\partial _{\nu } u\, (\m .\nabla u)\, d\sigma 
-2\int_{\Omega}\nabla u.\nabla(\m .\nabla u)\, d\x \, .
\end{equation*}
So, observing that 
$$\nabla u.\nabla(\m .\nabla u)
=\frac{1}{2}\m .\nabla (|\nabla u|^{2})+d|\nabla u|^{2}+(A\nabla u).\nabla u
\, ,$$ 
and since $A$ is skew-symmetric, we get 
\begin{equation*}
2\int_{\Omega}\Delta u\, (\m .\nabla u)\, d\x 
=\int_{\partial \Omega } 2\partial_{\nu } u\, (\m .\nabla u)\, d\sigma 
-2d\int_{\Omega }|\nabla u|^{2}\, d\x -\int_{\Omega }\m .\nabla
(|\nabla u|^{2})\, d\x \, .
\end{equation*}
With another use of Green-Riemann formula, we obtain the required result
for $\text{div} (\m )=nd$.
\end{proof}

\noindent We will now try to extend this result to the case of an element $u$ belonging less regular when $\Omega $ is smooth enough.

\subsection{Bi-dimensional case}
We begin by the plane case. It is the simplest case from the point of view of singularity theory and its understanding dates from Shamir \cite{Sh}. 
\begin{theorem}\label{T3}
Assume $n=2$. Under geometrical conditions \eqref{R} and \eqref{G'}, let 
$u\in \mathrm{H}^1(\Omega )$ such that
\begin{equation}\label{H}
\Delta u \in \mathrm{L}^2 (\Omega)\, , \quad
u_{/\partial \Omega_D } \in \mathrm{H}^{3/2}(\partial \Omega_D )\, , \quad 
\partial_{\nu} u_{/\partial \Omega_N } \in 
\mathrm{H}^{1/2} (\partial \Omega_N )\, .
\end{equation}
Then $\, 2\partial _{\nu }u\, (\m .\nabla u)-\m .\bn \, |\nabla u|^{2} $ belongs to
$\mathrm{L}^{1}(\partial \Omega )$ and there exist some coefficients 
$(c_{\x})_{\x \in \Gamma }$ such that
\begin{equation*}
2\int_{\Omega }\Delta u(\m .\nabla u)\, d\x =\int_{\partial \Omega }(2\partial
_{\nu } u\, (\m .\nabla u)-\m .\bn \, |\nabla u|^{2})\, d\sigma +\frac{\pi }{4}%
\sum_{\x \in \Gamma/\varpi_\x =\pi } c_{\x}^{2}\,(\m .\bt )(\x )\, .
\end{equation*}
\end{theorem}

\begin{proof}
We first begin by some considerations which will be used in the
general case too. It is a classical result that $u\in \mathrm{H}^{2}(\omega )$ 
for every open domain $\omega $ such that 
$\omega \Subset \overline{\Omega }\setminus \Gamma $. 
For sake of completeness, let us recall the proof.

\noindent A trace result shows that there exists 
$u_{R}\in \mathrm{H}^{2}(\omega )$ such that $u_{R}=u$ on $\partial \Omega_D $ 
and $\partial_{\nu }u_{R}=\partial_{\nu }u$ on $\partial \Omega_N $. 
Hence, setting $f=\Delta u_{R}-\Delta u\in\mathrm{L}^{2}(\Omega )$, 
$u_{S}=u-u_{R}$ satisfies
\begin{equation}\label{E}
\left\{ 
\begin{array}{l}
-\Delta u_{S}=f \\ 
u_{S}=0 \\ 
\partial _{\nu }u_{S}=0
\end{array}
\right. \left. 
\begin{array}{l}
\text{in }\Omega \, , \\ 
\text{on }\partial \Omega _{D} \, , \\ 
\text{on }\partial \Omega _{N} \, .
\end{array}
\right.
\end{equation}
$\bullet$ Now, if 
$\omega \Subset \Omega \setminus \Gamma \cup \partial \Omega_{D} $ and $\xi $ 
is a cut-off function such that $\xi =1$ on $\omega $ and
supp$(\xi )\subset \Omega $, then for a suitable 
$g\in\mathrm{L}^{2}(\Omega )$, $u_{\omega }=u_{S}\xi $ satisfies the Dirichlet 
problem
\begin{equation*}
\left\{ 
\begin{array}{l}
\Delta u_{\omega }=g \\ 
u_{\omega }=0
\end{array}
\right. \left. 
\begin{array}{l}
\text{on }\Omega \, , \\ 
\text{on }\partial \Omega \, ,
\end{array}
\right.
\end{equation*}
and using classical method of difference quotients (\cite{Gr}), one can
now conclude that $u_{\omega }\in \mathrm{H}^{2}(\Omega )$, hence 
$u_{S}\in \mathrm{H}^{2}(\omega )$.\\
$\bullet$ Else, if 
$\omega \Subset \Omega \setminus \Gamma \cup \partial \Omega_{N}$, 
and $\xi $ is a cut-off function such that $\xi =1$ on $\omega $ and 
supp$(\xi )\subset \Omega $, then for a suitable $g\in \mathrm{L}^{2}(\Omega )$,
$u_{\omega }=u_{S}\xi $ satisfies the Neumann problem
\begin{equation*}
\left\{ 
\begin{array}{l}
-\Delta u_{\omega }+u_{\omega }=g \\ 
\partial _{\nu }u_{\omega }=0
\end{array}
\right. \left. 
\begin{array}{l}
\text{on }\Omega \, , \\ 
\text{on }\partial \Omega \, ,
\end{array}
\right.
\end{equation*}
and, using similar argument, one gets $u_{S}\in \mathrm{H}^{2}(\omega )$.
\medskip

\noindent Let 
$\Omega _{\varepsilon }=\{\x\in \Omega \, /\, d(\x,\Gamma )>\varepsilon \}$.\\
By compactness of $\Omega _{\varepsilon }$, we get 
$u\in \mathrm{H}^{2}(\Omega_{\varepsilon })$. An application of Proposition
\ref{RR} to our particular situation gives us the following relation
\begin{equation*}
2\int_{\Omega _{\varepsilon }}\Delta u\, (\m .\nabla u)\, d\x =\int_{\partial
\Omega _{\varepsilon }}(2\partial _{\nu }u\, (\m .\nabla u)-\m .\bn \, |\nabla
u|^{2})\, d\sigma \, ,
\end{equation*}
and we will try to let $\varepsilon \rightarrow 0$. Using derivatives with
respect to $\bn $ and $\bt $, we get
\begin{equation*}
2\int_{\Omega _{\varepsilon }}\Delta u\, (\m .\nabla u)\, d\x =\int_{\partial
\Omega _{\varepsilon }}\m .\bn \, \left( (\partial _{\nu }u)^{2}-(\partial
_{\tau }u)^{2}\right)\, d\sigma +2\int_{\partial \Omega _{\varepsilon
}}\m .\bt \, (\partial _{\nu }u)(\partial _{\tau }u)\, d\sigma \, .
\end{equation*}%
First, since $\Delta u\in \mathrm{L}^{2}(\Omega )$ and 
$u\in \mathrm{H}^{1}(\Omega )$,
Lebesgue dominated convergence theorem immediately gives 
\begin{equation*}
\lim_{\varepsilon \rightarrow 0}\int_{\Omega _{\varepsilon }}\Delta
u\, (\m .\nabla u)\, d\x =\int_{\Omega }\Delta u\, (\m .\nabla u)\, d\x \, .
\end{equation*}
Now, we work on boundary terms. Let us introduce the following partition of 
$\partial \Omega _{\varepsilon }$: 
$\widetilde{\partial \Omega _{\varepsilon}}=
\partial \Omega _{\varepsilon }\cap \partial \Omega $, 
$\partial \Omega_{\varepsilon }^{\ast }=
\partial \Omega _{\varepsilon }\cap \Omega \ $ and use a
decomposition result due to Banasiak and Roach \cite{BR}: every variational 
solution of \eqref{E} can be split as a sum of singular functions.
There exist some coefficients $(c_{\x})_{\x\in \Gamma }$ and 
$u_{R}\in \mathrm{H}^{2}(\Omega )$ such that
\begin{equation}\label{D}
u=u_{R}+\sum_{\x\in \Gamma } c_{\x}U_{S}^{\x}=:u_{R}+u_{S}
\end{equation}
where $U_{S}^{\x}$ are singular functions which, in some
neighborhood of $\x\in \Gamma$, are defined in local polar coordinates 
(see Fig. 1) by
\begin{equation*}
U_{S}^{\x}(r,\theta )=\rho(r)\, r^{\frac{\pi }{2\varpi _{\x}}}\sin \bigl( 
\frac{\pi }{2\varpi _{\x}}\theta \bigr) \, .
\end{equation*}
with $\rho$ some cut-off function.\\
Using the density of $\mathcal{C}^{1}(\overline{\Omega })$ in 
$\mathrm{H}^{2}(\Omega )$, we will be able to assume that 
$u_{R}\in \mathcal{C}^{1}(\overline{\Omega })$.\\
Let us look at boundary terms on $\widetilde{\partial \Omega_{\varepsilon }}$
first. We first claim that for some constant $C>0$, 
\begin{equation*}
|\m .\bn |\leq Cd(.,\Gamma )\, .
\end{equation*}
In fact, if $\x\in \Omega $ and $\x_{1}\in \Gamma $ which satisfies $|\x-$ $%
\x_{1}|=d(\x,\Gamma )$, one gets
 $$\m .\bn (\x )=\m (\x ).(\bn (\x )-\bn (\x_{1} ))+(\m (\x )-\m (\x_{1} )).\bn
(\x_{1} )\quad (\text{observing that } \m .\bn (\x_{1} )=0)\, .$$ 
Hence, using the fact that $\bn $ is a piecewise $C^{1} $ function
(see Fig. 2), we get
 $$|\m .\bn (\x )|\leq (\Vert m\Vert _{\infty }\Vert \bn^{\prime } \Vert _{\infty
}+1)\, d(\x,\Gamma )\, .$$
Now, working in local coordinates, one gets
\begin{equation*}
d(\x,\Gamma )\, |\nabla u|^{2}\in \mathrm{L}^{\infty } (\partial \Omega ) \, .
\end{equation*}
Hence Lebesgue theorem implies 
\begin{equation*}
\lim_{\varepsilon \rightarrow 0}\int_{\widetilde{\partial \Omega
_{\varepsilon }}} \m .\bn \, \left( (\partial _{\nu }u)^{2}-(\partial _{\tau
}u)^{2}\right) \, d\sigma =\int_{\partial \Omega } \m .\bn \, \left( (\partial
_{\nu }u)^{2}-(\partial _{\tau }u)^{2}\right) \, d\sigma \, .
\end{equation*}
On the other hand, assumptions \eqref{H} give 
$$ \partial_{\nu } u_{/\partial \Omega_N } \in
\mathrm{H}^{1/2}(\partial \Omega_N) \, ,\quad
\partial_{\tau } u_{/\partial \Omega_N } \in
\mathrm{H}^{-1/2}(\partial \Omega_N) \, ,\quad
\partial_{\nu } u_{/\partial \Omega_D }\in 
\mathrm{H}^{-1/2}(\partial \Omega_D) \, ,\quad
\partial_{\tau } u_{/\partial \Omega_D }\in 
\mathrm{H}^{1/2}(\partial \Omega_D) \, .$$
Hence we get
\begin{equation*}
\int_{\widetilde{\partial \Omega _{\varepsilon }}} \m .\bt \, (\partial_{\nu }u)
(\partial_{\tau } u)\, d\sigma \longrightarrow \int_{\partial \Omega } \m .\bt \,
(\partial_{\nu } u)(\partial_{\tau } u) \, d\sigma \, ,\quad \hbox{ as }
\varepsilon \to 0\, .
\end{equation*}
Now, we have to consider the boundary term on $\partial \Omega _\varepsilon^* $, 
$I_{\varepsilon} (\nabla u)$.\\
It is a quadratic form with respect to $\nabla u$ and using \eqref{D}, one can 
decompose it as follows,
$$I_{\varepsilon} (\nabla u_R)+2J_{\varepsilon} (\nabla u_R,\nabla u_S)
+I_{\varepsilon} (\nabla u_S)\, ,$$
where $J_{\varepsilon}$ is the corresponding bilinear form.\\
Concerning $I_{\varepsilon} (\nabla u_R)$, regularity of $\m $ gives the estimate
\begin{equation*}
  | I_{\varepsilon} (\nabla u_R)| \leq C 
  \int_{\partial \Omega_{\varepsilon}^{\ast}} | \nabla u_R |^2 \, d\sigma
\end{equation*}
This term is $O(\varepsilon )$ since $\nabla u_R $ is bounded on $\Omega $.\\
For the term $I_{\varepsilon} (\nabla u_S)$, we first observe that, adjusting 
the cut-off functions, the supports of $u_{S}^{\x }$ and $u_{S}^{\y }$ are 
disjoint, provided that $\x \not =\y $.
Hence, using decomposition \eqref{D}, we can write
\begin{equation*}
I_{\varepsilon} (\nabla u_S)=\sum_{\x\in \Gamma} c_\x^2 
\int_{C_{\varepsilon }(\x )}(2\partial_{\nu } u^\x_S (\m .\nabla u^\x_S )
-\m .\bn \, |\nabla
u^\x_S|^{2} ) \, d\sigma\, .
\end{equation*}
If $\varpi _{\x}<\pi $, one gets 
\begin{equation*}
2\partial _{\nu }u_{S}^{\x} (\m .\nabla u_{S}^{\x} )- \m .\bn \, |\nabla
u_{S}^{\x}|^{2} =O(\varepsilon^{\frac{\pi }{\varpi _{\x}}-2})\, , \quad
\text{ on  } C_{\varepsilon }(\x)\, .
\end{equation*}
Hence, after integrating on $C_{\varepsilon }(\x )$, we get 
$\, \displaystyle{\lim_{\varepsilon \to 0} I_{1}^{\x}(\varepsilon  )= 0}$.

\begin{figure}[htb]
\centerline{
\begin{picture}(0,0)
\includegraphics{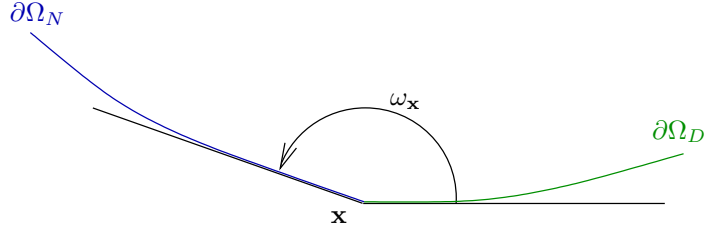}
\end{picture}
\setlength{\unitlength}{4144sp}
\begin{picture}(4061,1657)(1443,-2535)
\put(3300,-2300){\makebox(0,0)[lb]
{\smash{{\color[rgb]{0,0,0}$\x $}}}}
\put(5200,-1800){\makebox(0,0)[lb]
{\smash{{\color[rgb]{0,.56,0}$\partial \Omega_D $}}}}
\put(3650,-1600){\makebox(0,0)[lb]
{\smash{{\color[rgb]{0,0,0}$\omega_\x $}}}}
\put(1400,-1100){\makebox(0,0)[lb]
{\smash{{\color[rgb]{0,0,.69}$\partial \Omega_N $}}}}
\end{picture}
}
\caption{Shape of the boundary near an angular point $\x $.}
\end{figure}

\noindent If $\varpi_{\x} =\pi $, we will need the following identity
\begin{equation*}
2\partial _{\nu }u_{S}^{\x}\, (\m .\nabla u_{S}^{\x})-\m .\bn \, |\nabla
u_{S}^{\x}|^{2}=\frac{1}{4\varepsilon }(\m .\bt )(\x)\, ,\quad
\text{  on  }C_{\varepsilon }(\x)\, .
\end{equation*}
One can observe that $C_{\varepsilon }(\x)$ behaves as a half-circle when $
\varepsilon \rightarrow 0$. An integration gives
\begin{equation*}
\lim_{\varepsilon \to 0}\int_{C_{\varepsilon }(\x )}\bigl(2(\nu .\nabla
u_{S}^{\x })(\m .\nabla u_{S}^{\x})-\m .\bn \, |\nabla u_{S}^{\x}|^{2}\bigr)
\, d\sigma =\frac{\pi }{4} (\m .\bt )(\x).
\end{equation*}

\begin{figure}[htb]
\centerline{
\begin{picture}(0,0)
\includegraphics{sing1.pstex}
\end{picture}
\setlength{\unitlength}{4144sp}
\begin{picture}(3877,2726)(768,-2758)
\put(2450,-2600){\makebox(0,0)[lb]
{\smash{{\color[rgb]{0,0,0}$\bn (\x )$}}}}
\put(3700,-1450){\makebox(0,0)[lb]
{\smash{{\color[rgb]{0,0,0}$\bt (\x )$}}}}
\put(4500,-550){\makebox(0,0)[lb]
{\smash{{\color[rgb]{0,.56,0}$\partial \Omega_D $}}}}
\put(750,-600){\makebox(0,0)[lb]
{\smash{{\color[rgb]{0,0,.69}$\partial \Omega_N $}}}}
\put(2750,-1150){\makebox(0,0)[lb]
{\smash{{\color[rgb]{0,0,0}$\theta $}}}}
\put(2900,-700){\makebox(0,0)[lb]
{\smash{{\color[rgb]{1,0,0}$C_\varepsilon(\x )$}}}}
\put(3000,-150){\makebox(0,0)[lb]
{\smash{{\color[rgb]{0,0,0}$\bt (\y )$}}}}
\put(2750,-1750){\makebox(0,0)[lb]
{\smash{{\color[rgb]{0,0,0}$\bn (\y )$}}}}
\put(2350,-1550){\makebox(0,0)[lb]
{\smash{{\color[rgb]{0,0,0}$\x $}}}}
\put(2000,-740){\makebox(0,0)[lb]
{\smash{{\color[rgb]{0,0,0}$\y $}}}}
\end{picture}
}
\caption{Unit vectors $\bn (\x )$, $\bt (\x )$, $\bn (\y )$ and
$\bt (\y )$ when $\partial \Omega $ is regular at $\x $.}
\end{figure}

\noindent Finally, the bilinear term $J_{\varepsilon} (\nabla u_R, \nabla u_S) $
can be written entirely
$$\int_{\partial \Omega_{\varepsilon }^{\ast }} 
\partial_{\nu } u_R \, (\m. \nabla u_S) \, d\sigma 
+\int_{\partial \Omega_{\varepsilon }^{\ast }}
\partial_{\nu } u_S \, (\m. \nabla u_R) \, d\sigma 
-\int_{\partial \Omega_{\varepsilon }^{\ast }} 
(\m. \bn )\, (\nabla u_R . \nabla u_S )\, d\sigma \, . $$
Using the regularity of $\m $ and Cauchy-Schwarz inequality, we get an estimate
of the form
\begin{equation*}
  |J_{\varepsilon} (\nabla u_R, \nabla u_S) | \leq C \Bigl(
  \int_{\partial \Omega_{\varepsilon}^{\ast}} | \nabla u_R |^2 \, d\sigma
 \Bigr)^{1/2} \Bigl( \int_{\partial \Omega_{\varepsilon}^{\ast}} |
  \nabla u_S |^2 \, d\sigma \Bigr)^{1/2} \, .
\end{equation*}
We have seen that the first term in this inequality vanishes when 
$\varepsilon\rightarrow 0$.
For the second one, we now observe that, if $\varepsilon $ is small enough
$$\partial \Omega_{\varepsilon}^{\ast} = 
\bigsqcup_{\x \in  \Gamma} C_{\varepsilon} (\x)\, ,$$ 
where $C_{\varepsilon} (\x )$ is an arc of circle of radius $\varepsilon $
centered at $\x $. 
Then, we may write 
$$\int_{\partial \Omega_{\varepsilon}^{\ast}} | \nabla u_S |^2 \, d\sigma \leq 
2\sum_{\x ,\y \in \Gamma } c_\y^2 \int_{C_{\varepsilon }(\x )}|\nabla U_{S}^{\y}|^{2}
\, d\sigma\, .$$ 
A similar computation shows that, for $\x\in\Gamma$,
$\, \displaystyle{\int_{C_{\varepsilon }(\x)}|\nabla U_{S}^{\x}|^{2}\, d\sigma
=O(1)}$. Moreover, if $\x\neq \y $, $U_{S}^{\y }$ is bounded near $\x $, we get
$\, \displaystyle{\int_{C_{\varepsilon }(\x)} |\nabla U_{S}^{\y }|^{2}\, d\sigma
=O(\varepsilon)}$.
This completes the proof.
\end{proof}

\begin{remark}
The assumption $\mathcal{H}^{1}(\partial \Omega_{D})>0$ is not necessary in the above proof. We will now see why we need this assumption on the Dirichlet part in higher dimension.
\end{remark}

\subsection{General case}
We now state the result in general dimension.
\begin{theorem}\label{T4}
Assume $n\geq 3$. Under geometrical conditions \eqref{G} and \eqref{R}, let 
$u\in \mathrm{H}^1(\Omega )$ such that 
\begin{equation}\label{H'}
\Delta u \in \mathrm{L}^2 (\Omega) \, , \quad 
u_{/\partial \Omega_D } \in \mathrm{H}^{3/2}(\partial \Omega_D ) \, , \quad 
\partial_\nu u_{/\partial \Omega_N } \in \mathrm{H}^{1/2} (\partial \Omega_N )\, .
\end{equation}
Then, 
$\, 2\partial _{\nu }u\, (\m .\nabla u)-\m .\bn \, |\nabla u|^{2}$ belongs to
$\mathrm{L}^{1}(\partial \Omega )$ and there exists 
$\zeta \in \mathrm{H}^{1/2}(\Gamma )$ such that 
\begin{equation*}
2\int_{\Omega }\Delta u\, (\m .\nabla u)\, d\x 
=d(n-2)\int_{\Omega }|\nabla u|^{2}\, d\x 
+\int_{\partial \Omega } (2\partial _{\nu } u\, (\m .\nabla u)-\m .\bn \,
|\nabla u|^{2})\, d\sigma \\
+\int_{\Gamma }\m .\bt \, |\zeta |^{2} \, d\gamma \, .
\end{equation*}
\end{theorem}

\begin{proof}
We will essentially follow \cite{BLM}. As in the plane case, we set 
$\Omega _{\varepsilon }=\{\x\in \Omega ;$ $d(\x,\Gamma )>\varepsilon \}$. 
For any given $\varepsilon>0$, we may apply the identity of Proposition 
\ref{RR}
\begin{equation*}
2\int_{\Omega _{\varepsilon }} \Delta u\, (\m .\nabla u)\, d\x 
=d(n-2)\int_{\Omega_{\varepsilon }}|\nabla u|^{2}\, d\x 
+\int_{\partial \Omega_{\varepsilon }} (2\partial _{\nu }u\, (\m .\nabla u)
-\m .\bn \, |\nabla u|^{2})\, d\sigma \, ,
\end{equation*}
and we will again analyze the behavior of each term as
$\varepsilon\rightarrow 0$.

\noindent $\bullet $ First, since $\Delta u \in \mathrm{L}^{2}(\Omega )$ and 
$u\in \mathrm{H}^{1}(\Omega )$,
Lebesgue dominated convergence theorem immediately gives 
\begin{equation*}
\lim_{\varepsilon \to 0}\int_{\Omega_{\varepsilon }}\Delta u\, (\m .\nabla u)\, d\x 
=\int_{\Omega }\Delta u\, (\m .\nabla u)\, d\x \, , \quad
\lim_{\varepsilon \to 0}\int_{\Omega _{\varepsilon }}|\nabla u|^{2}\, d\x 
=\int_{\Omega } |\nabla u|^{2}\, d\x \, .
\end{equation*}
Below we shall consider boundary terms. We define 
$\, \widetilde{\partial\Omega _{\varepsilon }}
=\partial \Omega _{\varepsilon }\cap \partial \Omega $
and 
$\, \partial \Omega _{\varepsilon }^{\ast }
=\partial \Omega _{\varepsilon}\cap \Omega $ (see Fig. 3).

\begin{figure}[htb]
\centerline{
\begin{picture}(0,0)
\includegraphics{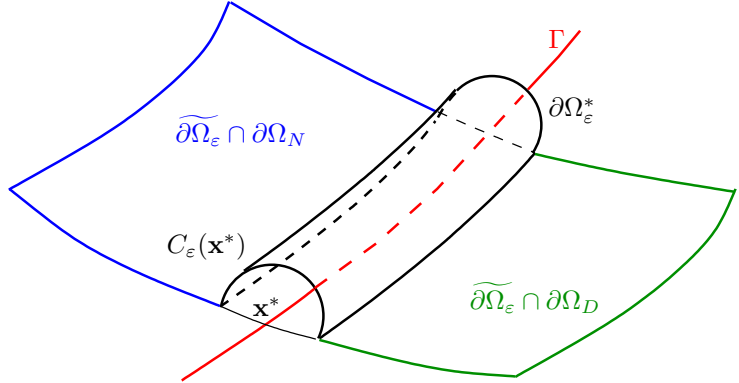}%
\end{picture}%
\setlength{\unitlength}{4144sp}%
\begingroup\makeatletter\ifx\SetFigFontNFSS\undefined%(dI+A)
\gdef\SetFigFontNFSS#1#2#3#4#5{%
  \reset@font\fontsize{#1}{#2pt}%
  \fontfamily{#3}\fontseries{#4}\fontshape{#5}%
  \selectfont}%
\fi\endgroup%
\begin{picture}(4345,2305)(988,-2408)
\put(4150,-811){\makebox(0,0)[lb]{\smash{
{\color[rgb]{0,0,0}$\partial\Omega_{\varepsilon}^{\ast}$}
}}}
\put(4150,-400){\makebox(0,0)[lb]{\smash{
{\color[rgb]{1,0,0}$\Gamma $}
}}}
\put(2400,-2000){\makebox(0,0)[lb]{\smash{
{\color[rgb]{0,0,0}$\x^{\ast}$}
}}}
\put(1891,-1625){\makebox(0,0)[lb]{\smash{
{\color[rgb]{0,0,0}$C_{\varepsilon}(\x^{\ast})$}
}}}
\put(1943,-962){\makebox(0,0)[lb]{\smash{
{\color[rgb]{0,0,1}$\widetilde{\partial \Omega_{\varepsilon}}\cap\partial \Omega_{N}$}
}}}
\put(3684,-1959){\makebox(0,0)[lb]{\smash{
{\color[rgb]{0,.56,0}$\widetilde{\partial \Omega_{\varepsilon}}\cap\partial \Omega_{D}$}
}}}
\end{picture}
}
\caption{Picture of $\partial\Omega_{\varepsilon}^{\ast}$ and $\widetilde{\partial \Omega_{\varepsilon}}$}
\end{figure}

\noindent $\bullet $ Let us consider boundary integral terms on 
$\widetilde{\partial\Omega _{\varepsilon }}$.\\
As well as in the plane case, there exists some constant $C>0$ such that 
$|\m .\bn |\leq C\, d(.,\Gamma )$. Thus, using the fact that 
$$ d(.,\Gamma )|\nabla u|^{2}\in \mathrm{L^{1}(\partial \Omega )}$$
(see \cite{BLM}, Proposition 3), we can use again Lebesgue theorem to conclude
that, as $\varepsilon \rightarrow 0$, 
\begin{equation*}
\int_{\widetilde{\partial \Omega _{\varepsilon }}}\m .\bn \, |\nabla
u|^{2}\, d\sigma \rightarrow \int_{\partial \Omega }\m .\bn \, |\nabla
u|^{2} \, d\sigma \, .
\end{equation*}
For the second integral, denoting by $\nabla_{\partial \Omega } $ the tangential
gradient along $\partial \Omega $, we write that 
\begin{equation*}
\partial_{\nu } u\, (\m .\nabla u) =\m .\bn \, |\partial_{\nu } u|^2
+\partial_{\nu } u\, (\m .\nabla_{\partial \Omega} u)\, .
\end{equation*}
The first term is integrable. The second one is, on $\partial \Omega_N$, the 
product of a $\mathrm{H^{1/2}}$ term by a $\mathrm{H^{-1/2}}$ one and, on 
$\partial \Omega_D$, the product of a $\mathrm{H^{-1/2}}$ term by a 
$\mathrm{H^{1/2}}$ one. Hence, Lebesgue theorem gives again, as 
$\varepsilon \rightarrow 0$, 
\begin{equation*}
\int_{\widetilde{\partial \Omega _{\varepsilon }}}\partial _{\nu }u\, (\m .\nabla
u)\, d\sigma \rightarrow \int_{\partial \Omega }\partial _{\nu }u\, (\m .\nabla
u)\, d\sigma \, .
\end{equation*}
$\bullet $ Let us now consider boundary integral terms on 
$\partial \Omega_{\varepsilon}^{\ast}$.\\
We assume that $\varepsilon \leqslant \varepsilon_0$ and we define
$\, \omega_{\varepsilon_0} : = \Omega \backslash \Omega_{\varepsilon_0}$.
As well as in the plane case, we can write
\begin{equation}\label{D'}
u = u_R + u_S
\end{equation}
where $u_S$ is the variational solution of some homogeneous mixed boundary 
problem and $u_R$ belongs to $\mathrm{H}^2 (\omega_{\varepsilon_0})$. Working by 
approximation if necessary, we
can suppose that $u_R \in \mathcal{C}^1 ( \overline{\omega_{\varepsilon_0}})$. 
Considering the same quadratic form as in the bi-dimensional case, this
leads to the following splitting 
$$\int_{\partial \Omega*_{\varepsilon
}}(2\partial _{\nu }u\, (\m .\nabla u)-\m .\bn \, |\nabla u|^{2})\, d\sigma =
I_{\varepsilon} (\nabla u_R) + I_{\varepsilon} (\nabla u_S) +
2 J_{\varepsilon} (\nabla u_R, \nabla u_S) \, . $$
Since $\nabla u_R \in \mathrm{L^{\infty} (\omega_{\varepsilon_0})}$ and 
$\mathcal{H}^{n -1} (\partial \Omega_{\varepsilon}^{\ast}) \rightarrow 0$ as 
$\varepsilon \rightarrow 0$, the first term 
$I_{\varepsilon} (\nabla u_R)$ clearly vanishes.

\noindent As above, the bilinear term $J_{\varepsilon} (\nabla u_R, \nabla u_S)$
is
$$\int_{\partial \Omega_{\varepsilon}^{\ast}} \partial_{\nu} u_R \, (\m. \nabla
u_S) \, d\sigma
+\int_{\partial \Omega_{\varepsilon}^{\ast}} \partial_{\nu} u_S \, (\m. \nabla
u_R) \, d\sigma 
-\int_{\partial \Omega_{\varepsilon}^{\ast}} (\m. \bn ) \, (\nabla u_R .\nabla
u_S) \, d \sigma \, .$$
Using the regularity of $\m $ and Cauchy-Schwarz inequality, we get an estimate
of the form
\begin{equation}\label{I}
  |J_{\varepsilon} (\nabla u_R, \nabla u_S) | \leqslant C \Bigl(
  \int_{\partial \Omega_{\varepsilon}^{\ast}} | \nabla u_R |^2 \, d\sigma
  \Bigr)^{1/2} \Bigl( \int_{\partial \Omega_{\varepsilon}^{\ast}} |
  \nabla u_S |^2 \, d\sigma \Bigr)^{1/2} \, .
\end{equation}
As above, it is clear that the first term vanishes as $\varepsilon
\rightarrow 0$.\\
 In order to analyze $I_{\varepsilon} (\nabla u_S)$ we will need further 
results.\\
To begin with, we introduce some notations.\\ 
Every $\x \in$ $\partial \Omega_{\varepsilon}^{\ast}$
belongs to a unique plane $\x^{\ast} + \left\langle \bt^{\ast}, \bn^{\ast}
\right\rangle$ (setting: $\bt^{\ast} = \bt (\x^{\ast})$, $\bn^{\ast} = \bn
(\x^{\ast})$) and more precisely to an arc-circle $C_{\varepsilon} (\x^{\ast})$
of center $\x^{\ast} \in \Gamma $ and of radius $\varepsilon $ 
(the figure is similar to Fig. 2 in the plane 
$\x^{\ast} + \left\langle \bt^{\ast}, \bn^{\ast}\right\rangle $ ). 
We define 
$$D_{\varepsilon} (\x^{\ast}) : = \omega_{\varepsilon} \cap (\x^{\ast} +
\left\langle \bt^{\ast}, \bn^{\ast} \right\rangle )\, .$$ 
For any $\x \in D_{\varepsilon_0} (\x^{\ast})$, we separate the derivatives of 
$u$ along the sub-manifold $\x - \x^{\ast} + \Gamma $ with the co-normal 
derivatives
\begin{equation}\label{DG}
\nabla u (\x )=\nabla_{\Gamma} u(\x )+\nabla_2 u(\x )\, ,\quad 
\nabla_{\Gamma} u(\x )\in T_{\x^{\ast}} \Gamma, \quad
\nabla_2 u(\x ) \in \left\langle \bt^{\ast}, \bn^{\ast} \right\rangle \, .
\end{equation}
Using methods of difference quotients (see for instance \cite{BLM}, Theorem 4),
one gets $\nabla_{\Gamma} u \in \mathrm{H^1 (\omega_{\varepsilon_0})}$ i.e.
$\nabla_{\Gamma} u_S \in \mathrm{H^1 (\omega_{\varepsilon_0})}$. We shall also
need the following result concerning the behavior of boundary integrals.

\begin{lemma}\label{LT}
Let $\varepsilon _{0}>0$. Assume that $u$ is such that 
$u=0$ on $\partial \omega_{\varepsilon _{0}}\cap \partial \Omega _{D}$,
$$ \forall \x^{\ast }\in \Gamma\, ,\quad u(\x^{\ast },.)\in
\mathrm{H}^{1}(D_{\varepsilon _{0}}(\x^{\ast })) \, ,$$ 
and 
$$\bigl( \x^{\ast }\longmapsto \Vert u(\x^{\ast },.)
\Vert_{\mathrm{H}^{1}(D_{\varepsilon_{0}(\x^{\ast }))}}\bigr) \in 
\mathrm{L}^{2}(\Gamma )\, .$$ 
Then there exists $C>0$
depending only on $\Omega $ such that, for any $\varepsilon $ sufficiently 
small,
\begin{equation*}
\int_{\Gamma }\Vert u(\x^{\ast },.)\Vert _{\mathrm{L^{2}(C_{\varepsilon }(\x^{\ast
}))}}^{2}\, d\gamma(\x^*) \leq 
C\varepsilon \int_{\Gamma }\Vert u(\x^{\ast },.)\Vert
_{\mathrm{H^{1}(D_{\varepsilon }(\x^{\ast }))}}^{2}\, d\gamma(\x^*) \, .
\end{equation*}
\end{lemma}

\noindent \textbf{Proof of Lemma \ref{LT}.}
We begin by changing coordinates as well as in \cite{BLM}. For every $%
\x_{0}^{\ast }\in \Gamma $, there exists $\rho _{0}>0$, a $\mathcal{C}^{2}$-
diffeomorphism $\Theta $ from an open neighborhood $W$ of 
$\x_{0}^{\ast }$ to 
$B(\rho_{0}):=B_{n-2}(\rho _{0})\times B_{2}(\rho _{0})$ (see Fig. 5) such that
\begin{eqnarray*}
\Theta (\x_{0}^{\ast }) &=&0\, , \\
\Theta (W\cap \Omega ) &=&\{ \y \in B(\rho _{0})\, / \, y_{n}>0\}\, , \\
\Theta (W\cap \partial \Omega _{D}) &=&\{ \y \in B(\rho
_{0})\, / \, y_{n-1}>0\, ,\, y_{n}=0\} \, , \\
\Theta (W\cap \partial \Omega _{N}) &=&\{ \y \in B(\rho
_{0})\, / \,y_{n-1}<0\, ,\, y_{n}=0\} \, , \\
\Theta (W\cap \Gamma ) &=&\{ \y \in B(\rho _{0})\, / \, y_{n-1}=0\, ,\, 
y_{n}=0\}:=\gamma
(\rho _{0}) \, .
\end{eqnarray*}

\begin{figure}[htb]
\centerline{
\begin{picture}(0,0)%
\includegraphics{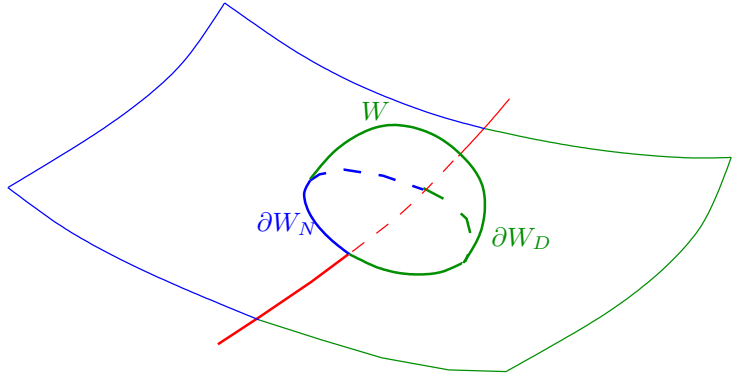}%
\end{picture}%
\setlength{\unitlength}{4144sp}%
\begingroup\makeatletter\ifx\SetFigFontNFSS\undefined%
\gdef\SetFigFontNFSS#1#2#3#4#5{%
  \reset@font\fontsize{#1}{#2pt}%
  \fontfamily{#3}\fontseries{#4}\fontshape{#5}%
  \selectfont}%
\fi\endgroup%
\begin{picture}(4305,2224)(997,-2314)
\put(3100,-800){\makebox(0,0)[lb]{\smash{{\color[rgb]{0,0.56,0}$W$}%
}}}
\put(3871,-1550){\makebox(0,0)[lb]{\smash{{\color[rgb]{0,.56,0}$\partial W_D$}%
}}}
\put(2483,-1464){\makebox(0,0)[lb]{\smash{{\color[rgb]{0,0,1}$\partial W_N$}%
}}}
\end{picture}
}
\caption{The set $W$.}
\end{figure}

\noindent Reducing $\varepsilon _{0}$ if necessary, we may assume that 
$D_{\varepsilon_{0}}(\x_{0}^*)\subset W$.\\
We then get, writing for $\x\in $ $W,$ $\
\Theta (\x)=(Y,\widetilde{y})\in \R^{n-2} \times \R^2 $ and 
$v:=u\circ \Theta ^{-1}$,
\begin{equation*}
\int_{W\cap \Gamma } \int_{C_{\varepsilon }(\x^{\ast })} u^{2} \, d\ell 
\, d\gamma (\x*)
=\int_{_{\gamma (\rho _0)}} \int_{\Theta (C_{\varepsilon }(\x^{\ast }))}
v^{2} \, d\ell(\widetilde{y}) \, dY \, .
\end{equation*}
Setting 
\begin{equation*}
B_{2}^{+}(\rho ):=\{\widetilde{y}=(y_{n-1},y_{n})\in B_{2}(\rho )\, /\, 
y_{n}>0\} \, ,\quad
C_{2}^{+}(\rho ):=\{\widetilde{y}=(y_{n-1},y_{n})\in \partial B_{2}(\rho )\, /\,
y_{n}>0\}\, ,
\end{equation*}
we first observe that we can choose $\rho _{\x^{\ast }}$ such that 
$\{ Y\}\times B_{2}^{+}(\rho )\subset \Theta (D_{\varepsilon }(\x^{\ast }))$.
Hence denoting by $\pi _{2}$ the projection on 
$\{0_{\mathbb{R}^{n-2}}\}\times\mathbb{R}^2$, the change of variables 
\begin{eqnarray*}
     \pi _{2}(\Theta (C_{\varepsilon }(\x^{\ast }))) &\longrightarrow & C_{2}^{+}(\rho )\\
     \widetilde{y} &\longmapsto& z=\rho \frac{\widetilde{y}}{|\widetilde{y} |}
\end{eqnarray*}
gives the estimate
\begin{equation}\label{IL}
\int_{\Theta (C_{\varepsilon }(\x^{\ast }))}v(Y,\widetilde{y})^2 \, 
d\ell(\widetilde{y})\leq
C\int_{C_{2}^{+}(\rho )}v(Y,z)^2\, d\ell(z)
\end{equation}
for a constant $C$ depending only on $\x_{0}^{\ast }$.

\noindent We will now estimate this latter integral in terms of $\Vert \nabla
_{2}v\Vert _{\mathrm{L}^{2}(\{y'\}\times B_{2}^{+}(\rho ))}$. Setting $v_{\rho
}(\widetilde{y}):=v(Y,\widetilde{y})$, one gets 
$\nabla v_{\rho }\in \mathrm{L}^{2}(B_{2}^{+}(1))$
and
\begin{equation*}
\Vert \nabla v_{\rho }\Vert _{\mathrm{L}^{2}(B_{2}^{+}(1))}
=\Vert \nabla _{2}v\Vert_{\mathrm{L}^{2}(\{y'\}\times B_{2}^{+}(\rho ))} \, ,\quad
\Vert v_{\rho }\Vert_{\mathrm{L}^{2}(C_{2}^{+}(1))}
=\rho^{-\frac{1}{2}}\Vert v\Vert _{\mathrm{L}^{2}(\{y'\}\times C_{2}^{+}(\rho ))} \, .
\end{equation*}
Observing that $v_{\rho }=0$ on $B_{2}^{++}(1):=\{(y_{n-1},y_{n})\in
B_{2}^{+}(1)\, /\, y_{n}>0\}$, trace theorem and Poincar\'{e} inequality give,
for some universal constant $C>0$, the estimate
\begin{equation*}
\int_{C_{2}^{+}(\rho )}v^{2}(y',\widetilde{y})\, d \ell(\widetilde{y})\leq C\rho
\Vert \nabla _{2}v\Vert _{\mathrm{L}^{2}(\{Y\}\times B_{2}^{+}(\rho ))}^{2} \, .
\end{equation*}
Hence, thanks to \eqref{IL}, one gets
\begin{equation*}
\int_{\Theta (C_{\varepsilon }(\x^{\ast }))}v^{2}(Y,\widetilde{y})\, 
d\ell(\widetilde{y})\leq C\rho_{\x^{\ast }} \Vert \nabla_{2} v
\Vert_{\mathrm{L}^{2}(\{ Y\} \times B_{2}^{+}(\rho_{\x^{\ast }}))}^{2}\, .
\end{equation*}

\begin{figure}[htb]
\centerline{
\begin{picture}(0,0)
\includegraphics{sing2.pstex}
\end{picture}
\setlength{\unitlength}{4144sp}
\begin{picture}(5877,1790)(892,-1846)
\put(1877,-1350){\makebox(0,0)[lb]
{\smash{{\color[rgb]{0,0,0}$\x^* $}}}}
\put(3021,-1200){\makebox(0,0)[lb]
{\smash{{\color[rgb]{0,.56,0}$\partial \Omega_D $}}}}
\put(5750,-1400){\makebox(0,0)[lb]
{\smash{{\color[rgb]{0,0,0}$O$}}}}
\put(2000,-800){\makebox(0,0)[lb]
{\smash{{\color[rgb]{1,0,0}$\Omega $}}}}
\put(4000,-160){\makebox(0,0)[lb]
{\smash{{\color[rgb]{0,0,0}$\Theta(\x^*,.) $}}}}
\put(907,-1200){\makebox(0,0)[lb]
{\smash{{\color[rgb]{0,0,.69}$\partial \Omega_N $}}}}
\put(1994,-1700){\makebox(0,0)[lb]
{\smash{{\color[rgb]{0,0,0}$\bn (\x^*) $}}}}
\end{picture}
}
\caption{The ${\cal C}^2$-diffeomorphism $\Theta(\x^*,.) $ in the plane 
$\x^* + \left\langle \bt^*, \bn^*\right\rangle$.}
\end{figure}

\noindent Observing that $\rho _{\x^{\ast }}$ is uniformly $O(\varepsilon )$ on 
$ W\cap \Gamma $ and the diffeomorphism $\Theta (\x^{\ast },.)$ (see Fig. 6),
we can conclude that, for some constant $C_{\x_{0}^{\ast }}$ depending only on
$\x_{0}^{\ast }$
\begin{equation*}
\int_{\Theta (C_{\varepsilon }(\x*))}v^{2}(Y,\widetilde{y})\, d\ell(\widetilde{y})
\leq
C_{\x_{0}^{\ast }}\varepsilon \Vert u(\x^{\ast },.)\Vert _{\mathrm{H}^{1}(\Theta
^{-1}(\{Y\}\times B_{2}^{+}(\rho )))}^{2} \leq 
C_{\x_{0}^{\ast }}\varepsilon \Vert u(\x^{\ast },.)\Vert
_{\mathrm{H}^{1}(D_{\varepsilon }(\x^{\ast }))}^{2} \, .
\end{equation*}
Hence, after an integration on $W\cap \Gamma $
\begin{equation*}
\int_{W\cap \Gamma } \int_{C_{\varepsilon }(\x^{\ast })}u^{2} \, d\ell 
\, d\gamma (\x^{\ast })\leq 
C_{\x_{0}^{\ast }} \varepsilon \int_{W\cap \Gamma }\Vert u(\x^{\ast},.)
\Vert _{\mathrm{H}^{1}(D_{\varepsilon }(\x^{\ast }))}^{2} \, d\gamma(\x^*) \, .
\end{equation*}
We finally complete the proof by using a partition of unity on the open sets 
$ (W_{\x_{0}^{\ast }})_{\x_{0}^{\ast }\in \Gamma }$.
\\ \textit{End of proof of Lemma \ref{LT}.} $\blacksquare$

\noindent Let us come back to our problem. Using \eqref{DG} for $u_S$, Pythagore
theorem gives
$$\int_{\partial \Omega_{\varepsilon}^{\ast}} | \nabla u_S |^2 \, d\sigma 
=\int_{\partial \Omega_{\varepsilon}^{\ast}} |\nabla_{\Gamma} u_S |^2 \, 
d\sigma 
+\int_{\partial \Omega_{\varepsilon}^{\ast}} | \nabla_2 u_S |^2 \, 
d\sigma \, . $$
Applying Lemma \ref{LT} to $\nabla_{\Gamma} u_S$, we get that the first term 
vanishes as $\varepsilon \rightarrow 0$. As well as in the bi-dimensional case,
we will
see that the second term above is bounded, using more information on $u_S$.

\noindent Thanks to \cite{BLM} (Theorem 4) and Borel-Lebesgue theorem, we may 
write
\begin{equation}\label{S}
 u_S (\x) = \eta (\x^{\ast}) U_S (\x - \x^{\ast}):= \eta (\x^{\ast})
  U^{\x^{\ast}}_S (\x) \, ,\quad \text{ on } \omega_{\varepsilon_0} \, ,
\end{equation}
with $U_S$ locally diffeomorphic to Shamir function, and $\eta
\in H^{1/2} (\Gamma)$. We then get, thanks to Fubini theorem
$$\int_{\partial \Omega_{\varepsilon}^{\ast}} | \nabla_2 u_S |^2 \, d\sigma 
=\int_{\Gamma} \eta (\x^{\ast})^2 \int_{C_{\varepsilon} (\x^{\ast})} 
|\nabla_2 U^{\x^{\ast}}_S |^2 \, d\ell \, d\gamma (\x^{\ast})\, ,$$
and, as well as in the bi-dimensional case, we show that this term is bounded 
by $O (1) \left\| \eta \right\|^2_{L^2 (\Gamma)}$. We have now proven that the 
second term in \eqref{I} is bounded, that is
$$J_{\varepsilon} (\nabla u_R) \rightarrow 0\, ,\quad \text{ as } \varepsilon 
\rightarrow 0\, .$$

\noindent To treat the last term $I_{\varepsilon} (\nabla u_S)$, we will use 
similar tools. The splitting \eqref{D'} for $u_S$ gives us
\[ I_{\varepsilon} (\nabla u_S) = I_{\varepsilon} (\nabla_2 u_S) +
   I_{\varepsilon} (\nabla_{\Gamma} u_S) + 2J_{\varepsilon} (\nabla_2 u_S,
   \nabla_{\Gamma} u_S) \, . \]
As above, the term $\displaystyle{I_{\varepsilon} (\nabla_{\Gamma} u_S)}$ is 
estimated by
$\, \displaystyle{\int_{\partial \Omega_{\varepsilon}^{\ast}} | \nabla_{\Gamma}
 u_S |^2\, d\sigma}$. It then vanishes for $\varepsilon \rightarrow 0$.\\ 
The bilinear term is estimated by
$$\Bigl( \int_{\partial \Omega_{\varepsilon}^{\ast}} | \nabla_2 u_S |^2 \, 
d\sigma \Bigr)^{1/2} 
\Bigl( \int_{\partial \Omega_{\varepsilon}^{\ast}} | \nabla_{\Gamma} u_S |^2 
\, d\sigma \Bigr)^{1/2} \, ,$$
it then tends to zero since the first term is bounded and the second one 
vanishes for $\varepsilon \rightarrow 0$.

\noindent For the last term $I_{\varepsilon} (\nabla_2 u_S)$, we use \eqref{S} 
and Fubini theorem to write it
$$\int_{\Gamma} \eta (\x^{\ast})^2 \int_{C_{\varepsilon} 
\bigl((\x^{\ast})} 2 (\bn . \nabla_2 U^{\x^{\ast}}_S) (\m. \nabla_2 U^{\x^{\ast}}_S)
-\m. \bn \, |\nabla_2 U^{\x^{\ast}}_S |^2 \bigr) \, d\ell \, 
d\gamma (\x^{\ast})\, .$$
We first work in the plane $\x^{\ast} + \left\langle
\bt^{\ast}, - \bn^{\ast} \right\rangle$ and, as above, we get
\[ \lim_{\varepsilon \to 0} \int_{C_{\varepsilon} (\x^{\ast})} (2 (\bn .
   \nabla_2 U^{\x^{\ast}}_S) (\m. \nabla_2 U^{\x^{\ast}}_S) -\m. \bn \, | \nabla_2
   U^{\x^{\ast}}_S |^2) \, d\ell = \frac{\pi}{4} \m (\x^{\ast}). \bt (\x^{\ast}) . \]
Moreover, for any $\varepsilon>0$, this integral term on 
$C_{\varepsilon} (\x^{\ast})$ is dominated by 
$\displaystyle{\frac{\pi}{2} \|\m\|_{\infty}\in L^1(\Gamma)}$.
So dominated convergence theorem applies and finally
\[ \lim_{\varepsilon \rightarrow 0} I_{\varepsilon} (\nabla_2 u_S) =
   \frac{\pi}{4} \int_{\Gamma} \eta^2 \m. \bt \, d\gamma \,  . \]
The proof is now complete with 
$\displaystyle{\zeta = \frac{\sqrt{\pi}}{2} \eta }$.
\end{proof}

\noindent We will now apply Rellich relation to the
stabilization of solutions of $(S)$.

\section{Proof of linear and non-linear stabilization}

We begin by writing the following consequence of Section 2.
\begin{corollary}\label{C}
Assume that $t\mapsto(u(t),u'(t))$ is a strong solution of $(S)$ and that the 
geometrical additional assumption \eqref{S1} if $n\geq3$ or \eqref{S2} if $n=2$
holds.
Then for every time $t$, $u(t)$ satisfies
\begin{equation*}
2\int_{\Omega }\Delta u (\m .\nabla u)\, d\x  \leq 
d(n-2)\int_{\Omega }|\nabla u|^{2}\, d\x 
+\int_{\partial \Omega }(2\partial _{\nu }u\, (\m .\nabla u)-\m .\bn \, |\nabla
u|^{2})\, d\sigma \, .
\end{equation*}
\end{corollary}
\begin{proof}
Indeed, under theses hypotheses, for each time $t$, 
$(u(t),u'(t))\in D(\mathcal{A})$ so that $u(t)$ satisfies  \eqref{H} or 
\eqref{H'}. The corollary is then an application of Theorem \ref{T3} or 
\ref{T4}.
\end{proof}

\noindent We will be able to prove Theorems \ref{T1} and \ref{T2} showing 
that, for $\displaystyle{\alpha=\frac{p-1}{2}}$, one can apply the following 
result \cite{Ko}.

\begin{proposition}\label{K}
Let $E:\mathbb{R}_{+}\rightarrow \mathbb{R}_{+}$ a non-increasing function such that there exists $\alpha \geq 0$ and $%
C>0$ which fulfills%
\begin{equation*}
\forall t\geq 0\text{, }\int_{t}^{\infty }E^{\alpha +1}(s)\, ds\leq CE(t).
\end{equation*}%
Then, setting $T=CE^{\alpha }(0)$, one gets
\begin{eqnarray*}
    \text{if } \alpha=0,  &\forall t \geq T,& \, E (t) \leq E (0) \exp
    \Bigl( 1 - \frac{t}{T}\Bigr) \,  ,\\
    \text{if } \alpha>0,  &\forall t \geq T,&  E (t) \leq E (0) \Bigl(
    \frac{T + \alpha T}{T + \alpha t} \Bigr)^{1/\alpha}.
  \end{eqnarray*}
\end{proposition}

\noindent We come back to our proof now.
\\
\\
\begin{proof}
Following \cite{Ko} and \cite{CR}, we will prove the estimates for 
$(u_{0},u_{1})\in D(\mathcal{A})$ which, using density of the domain, will be 
sufficient to get the result for all solutions.\\
Setting $Mu=2\m .\nabla u+d(n-1)u$, we prove the following result.

\begin{lemma}\label{IPP}
For any $0\leq S<T<\infty $, one gets
\begin{eqnarray*}
2d\int_{S}^{T} E^{\frac{p+1}{2}}\, dt 
&\leq 
&-\Bigl[ E^{\frac{p-1}{2} }\int_{\Omega }u^{\prime }Mu \, d\x \Bigr]_{S}^{T}
+\frac{p-1}{2}\int_{S}^{T}E^{\frac{p-3}{2}}E^{\prime }
\int_{\Omega }u^{\prime } Mu \, d\x \, dt \\
&
&+\int_{S}^{T} E^{\frac{p-1}{2}} \int_{\partial \Omega _{N}}\m .\bn \,
\bigl( (u^{\prime })^{2}-|\nabla u|^{2}-g(u^{\prime })Mu\bigr) \, d\sigma \, dt 
\, .
\end{eqnarray*}
\end{lemma}

\noindent \textbf{Proof of Lemma \ref{IPP}.}
Using the fact that $u$ satisfies $(S)$ and observing that $u^{\prime
\prime }Mu=(u^{\prime }Mu)^{\prime }-u^{\prime }Mu^{\prime }$, an
integration by parts gives
\begin{eqnarray*}
0 
&=
&\int_{S}^{T} E^{\frac{p-1}{2}} \int_{\Omega } (u^{\prime \prime }-\Delta u) Mu
\, d\x \, dt \\
&=&\Bigl[ E^{\frac{p-1}{2}} \int_{\Omega } u^{\prime } Mu\, d\x \Bigr]_{S}^{T}
-\frac{p-1}{2} \int_{S}^{T} E^{\frac{p-3}{2}} E^{\prime } 
\int_{\Omega} u^{\prime } Mu\, d\x \, dt 
-\int_{S}^{T} E^{\frac{p-1}{2}} \int_{\Omega }(u^{\prime } Mu^{\prime }+\Delta u Mu)
\, d\x \, dt\, .
\end{eqnarray*}
Corollary \ref{C} now gives
\begin{equation*}
\int_{\Omega }\Delta u\, Mu\, d\x  \, \leq \, 
d(n-1)\int_{\Omega } \Delta u\, u\, d\x 
+d(n-2)\int_{\Omega } |\nabla u|^{2} \, d\x 
+\int_{\partial \Omega } (2\partial _{\nu } u\, (\m .\nabla u)-\m .\bn \,
|\nabla u|^{2}) \, d\sigma \, .
\end{equation*}
hence, Green-Riemann formula leads to
\begin{equation*}
\int_{\Omega } \Delta u\, Mu \, d\x \leq 
-d\int_{\Omega } |\nabla u|^{2} \, d\x 
+\int_{\partial \Omega } (\partial _{\nu } u\, Mu-\m .\bn \, |\nabla
u|^{2}) \, d\sigma \, .
\end{equation*}
Using boundary conditions and the fact that $\nabla u=\partial _{\nu
}u\, \bn $ on $\partial \Omega _{D}$, we get
\begin{equation*}
\int_{\Omega } \Delta u\, Mu \, d\x \leq 
-d\int_{\Omega } |\nabla u|^{2} \, d\x 
-\int_{\partial \Omega _{N}} \m .\bn \, \bigl( g(u^{\prime })\, Mu+|\nabla
u|^{2} \bigr) \, d\sigma \, .
\end{equation*}
On the other hand, using $\text{div}(\m )=nd$, another use of Green formula
gives us
\begin{equation*}
\int_{\Omega } u^{\prime } \, Mu^{\prime } \, d\x 
=-d\int_{\Omega } |u^{\prime } |^{2} \, d\x
+\int_{\partial \Omega _{N}} \m .\bn \, |u^{\prime }|^{2}\, d\sigma \, .
\end{equation*}
\textit{End of proof of Lemma \ref{IPP}.} $\blacksquare$

\noindent Coming back to our problem, Young inequality gives
 $$\Bigl| \int_{\Omega }u^{\prime } \, Mu \, d\x \Bigr| \leq
CE(t)\, .$$
Lemma \ref{IPP} shows that
\begin{eqnarray*}
2d\int_{S}^{T} E^{\frac{p+1}{2}} \, dt 
&\leq 
&C\bigl( E^{\frac{p+1}{2}}(T)+E^{\frac{p+1 }{2}}(S)\bigr)
+C\int_{S}^{T}E^{\frac{p-1}{2}}E^{\prime } \, dt \\
&
&+\int_{S}^{T}E^{\frac{p-1}{2}} \int_{\partial \Omega _{N}} \m .\bn \,
\bigl( |u^{\prime }|^{2}-|\nabla u|^{2}-g(u^{\prime })Mu \bigr)\, d\sigma 
\, dt \, .
\end{eqnarray*}
For simplicity, let $d\sigma _{m} =\m .\bn \, d\sigma $. 
Observing that 
$\displaystyle{E^{\prime }(t)
=-\int_{\partial \Omega _{N}} g(u^{\prime }) u^{\prime } \, d\sigma_{m}\leq 0}$,
we get, for a constant $C>0$ independent of $E(0)$ if $p=1$,
\begin{equation*}
2d\int_{S}^{T}E^{\frac{p+1}{2}} \, dt\leq
CE(S)+\int_{S}^{T}E^{\frac{p-1}{2} } \int_{\partial \Omega _{N}} \bigl(
|u^{\prime }|^{2}-|\nabla u|^{2}-g(u^{\prime })Mu \bigr) \, d\sigma_{m} \, dt\,  .
\end{equation*}
Using the definition of $Mu$ and Young inequality, we get, for any 
$\varepsilon >0$,
\begin{equation*}
2d\int_{S}^{T}E^{\frac{p+1}{2}}\, dt\leq 
CE(S)+\int_{S}^{T}E^{\frac{p-1}{2} } \int_{\partial \Omega _{N}}
\Bigl( [u^{\prime }|^{2}+\Bigl( \Vert \m \Vert_{\infty }^{2}
+\frac{d^{2}(n-1)^{2}}{4\varepsilon } \Bigr) g(u^{\prime })^{2}
+\varepsilon u^{2} \Bigr) \, d\sigma_{m} \, dt\, .
\end{equation*}
Now, using Poincar\'{e} inequality, we can choose $\varepsilon >0$ such
that
\begin{equation*}
\varepsilon \int_{\partial \Omega_{N} }\m .\bn \, u^{2}\, d\sigma \leq 
\frac{d}{2} \int_{\Omega } |\nabla u|^{2}\, d\x \leq dE \, . 
\end{equation*}
So we conclude
\begin{equation*}
d\int_{S}^{T}E^{\frac{p+1}{2}}\, dt\leq 
CE(S)+C\int_{S}^{T}E^{\frac{p-1}{2} } \int_{\partial \Omega _{N}}\bigl(
(u^{\prime })^{2}+g(u^{\prime })^{2}\bigr) \, d\sigma_{m} \, dt\, .
\end{equation*}
We split $\partial \Omega _{N}$ to bound the last term of the above estimate
 $$\partial \Omega _{N}^{1}=\{\x\in
\partial \Omega _{N};\, |u^{\prime }(\x)|>1\}, \quad \partial \Omega
_{N}^{2}=\{\x\in \partial \Omega _{N}; \, |u^{\prime }(\x)|\leq 1\} \, .$$

\noindent Using \eqref{F1} and \eqref{F2}, we get
\begin{equation*}
\int_{S}^{T}E^{\frac{p-1}{2}} \int_{\partial \Omega
_{N}^{1}}\left( |u^{\prime }|^{2}+g(u^{\prime })^{2}\right) \, d\sigma _{m}
\, dt \leq
C\int_{S}^{T}E^{\frac{p-1}{2}} \int_{\partial \Omega _{N}}u^{\prime}
g(u^{\prime })\, d\sigma_{m} \, dt \leq CE(S)\, ,
\end{equation*}
where $C$ neither depend on $E(0)$ if $p=1$.

\noindent On the other hand, using \eqref{F1}, \eqref{F2}; Jensen inequality
 and boundedness of $\m $, one
successively obtains
\begin{equation*}
\int_{\partial \Omega _{N}^{2}}\bigl( (u^{\prime })^{2}+g(u^{\prime })^{2}
\bigr) \, d\sigma_{m} \leq
C\int_{\partial \Omega _{N}^{2}}(u^{\prime }g(u^{\prime}))^{2/(p+1)}\, 
d\sigma_{m} \leq 
C\Bigl( \int_{\partial \Omega _{N}^{2}}
u^{\prime }g(u^{\prime })\, d\sigma _{m}\Bigr)^{\frac{2}{p+1}} \leq 
C(-E^{\prime })^{\frac{2}{p+1}} \, .
\end{equation*}
Hence, using Young inequality again, we get for every $\varepsilon >0$%
\begin{equation*}
\int_{S}^{T}E^{\frac{p-1}{2}}\int_{\partial \Omega_{N}^{2}}
\bigl( (u^{\prime })^{2}+g(u^{\prime })^{2} \bigr) \, d\sigma _{m} \, dt \leq
\int_{S}^{T}(\varepsilon E^{\frac{p+1}{2}}-C(\varepsilon )E^{\prime })\, dt
\leq \varepsilon \int_{S}^{T}E^{\frac{p+1}{2}}\, dt+C(\varepsilon )E(S)\, .
\end{equation*}
Finally, we get, for some $C(\varepsilon )$ and $C$ independent of $E(0)$
if $p=1$
\begin{equation*}
d\int_{S}^{T}E^{\frac{p+1}{2}}dt\leq C(\varepsilon )E(S)+\varepsilon
C\int_{S}^{T}E^{\frac{p+1}{2}}dt \, .
\end{equation*}
Choosing now $\displaystyle{\varepsilon C \leq \frac{d}{2}}$, Theorems 1 and 2
result from Proposition 8.
\end{proof}

\section{Examples and numerical results}

\subsection{Examples}
We here consider the case when $\Omega $ is a plane convex polygonal domain.
The normal unit vector pointing outward of $\Omega $ is piecewise
constant and the nature of boundary conditions involved by the multiplier
method can be determined on each edge, independently of other edges.\\
Along each edge, vector $\bn $ is constant and the boundary conditions are
defined by the sign of
$$ \m (\x ).\bn (\x )=(R_\theta (\x -\x_0 )).\bn (\x )
=(\x -\x_0 ).R_{-\theta} (\bn (\x ))\, .$$
Hence we build $\bn $, $R_{-\theta} (\bn ) $ and we we can determine the sign of
above coefficient with respect to the position of $\x_0 $. To this end, we
construct two straight lines, orthogonal with respect to $R_{-\theta} (\bn )$ so
that each of them contains one vertex of the considered edge.\\
This determines a belt and if $\x_0 $ belongs to this belt, we obtained mixed
boundary conditions along this edge, if $\x_0 $ does not belong to this belt,
then we get Dirichlet or Neumann boundary conditions along whole the edge
(see Figure 6).
\begin{figure}[htb]
\centerline{
\begin{picture}(0,0)
\includegraphics{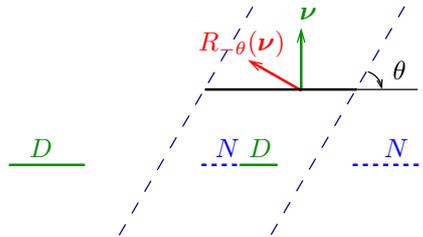}
\end{picture}
\setlength{\unitlength}{4144sp}
\begingroup\makeatletter\ifx\SetFigFont\undefined
\gdef\SetFigFont#1#2#3#4#5{
  \reset@font\fontsize{#1}{#2pt}
  \fontfamily{#3}\fontseries{#4}\fontshape{#5}
  \selectfont}
\fi\endgroup
\begin{picture}(2529,1426)(1069,-949)
\put(2600,354){\makebox(0,0)[lb]{\smash{
{\color[rgb]{0,.56,0}$\bn $}
}}}
\put(3150, -1){\makebox(0,0)[lb]{\smash{
{\color[rgb]{0,0,0}$\theta $}
}}}
\put(2000,169){\makebox(0,0)[lb]{\smash{
{\color[rgb]{1,0,0}$R_{-\theta } (\bn )$}
}}}
\put(3100,-466){\makebox(0,0)[lb]{\smash{
{\color[rgb]{0,0,1}$N$}
}}}
\put(2100,-466){\makebox(0,0)[lb]{\smash{
{\color[rgb]{0,0,1}$N$}
}}}
\put(2300,-466){\makebox(0,0)[lb]{\smash{
{\color[rgb]{0,.56,0}$D$}
}}}
\put(1000,-466){\makebox(0,0)[lb]{\smash{
{\color[rgb]{0,.56,0}$D$}
}}}
\end{picture}
}
\caption{Boundary conditions along some edge depending on the position of $\x_0 $.}
\end{figure}
\\
Performing this method for a square, $\Omega =(0,1)^2 $, we show in Figure 7 the
different cases of boundary conditions depending on the position of $\x_0 $.
Three main cases are considered
\begin{enumerate}
\item $\displaystyle{0<\theta <\frac{\pi}{4} }\, :$ above belts controlling
opposite edges have a non-empty intersection, which is a belt of positive
thickness,
\item $\displaystyle{\theta =\frac{\pi}{4} }\, :$ this intersection is a 
straight line,
\item $\displaystyle{\frac{\pi}{4} <\theta <\frac{\pi}{2}}\, :$ the 
intersection is empty.
\end{enumerate}
The case when $\theta $ is negative can be easily deduced by symmetry.\\
In the three above cases, there are four angular sectors (shaded areas in
Figure 7) such that if $\x_0 $ belongs to one of them, then geometrical
condition \eqref{S2} is satisfied.
\begin{figure}[htb]
\centerline{
\epsfig{figure=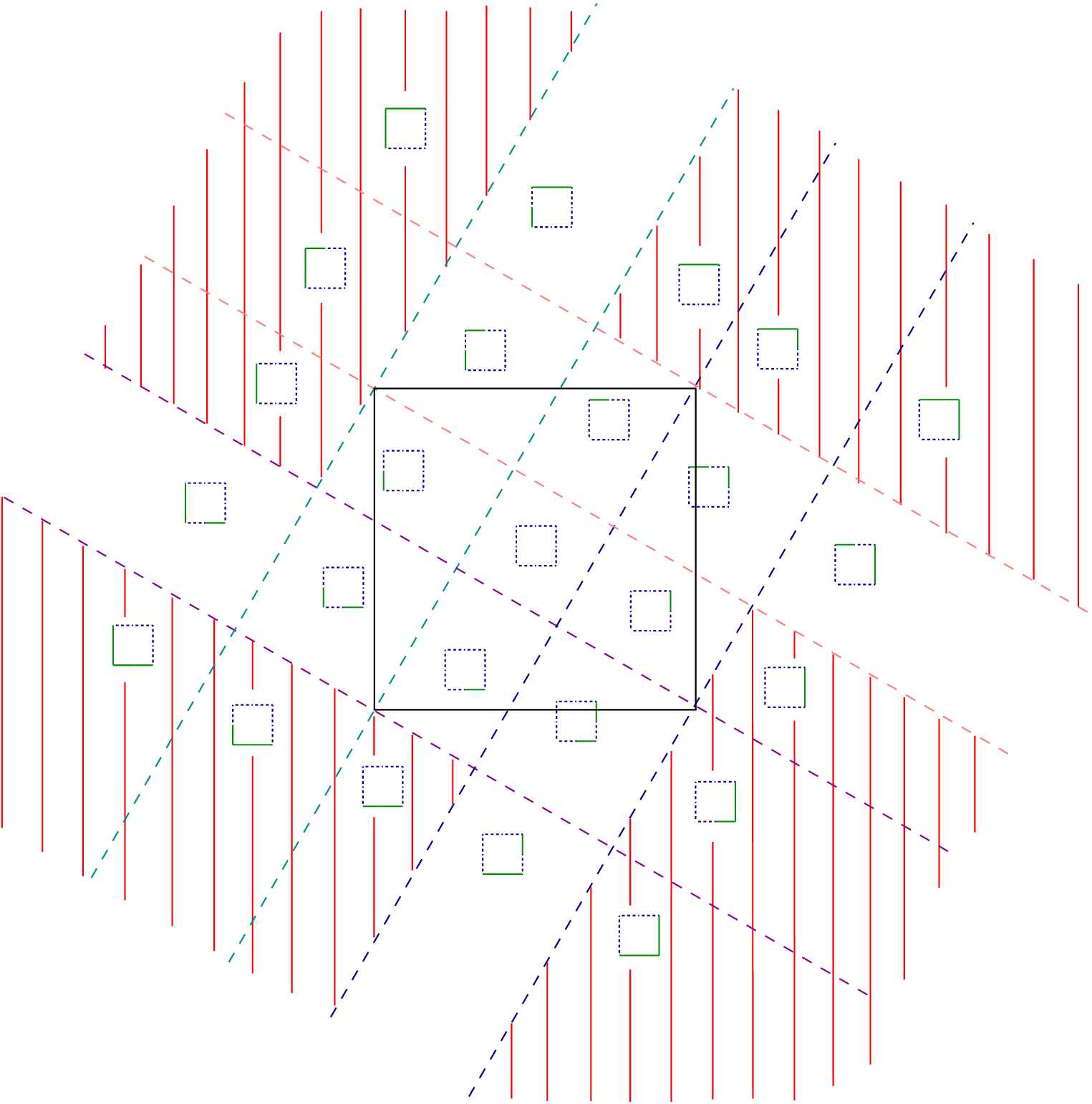,scale=0.18}
\qquad
\epsfig{figure=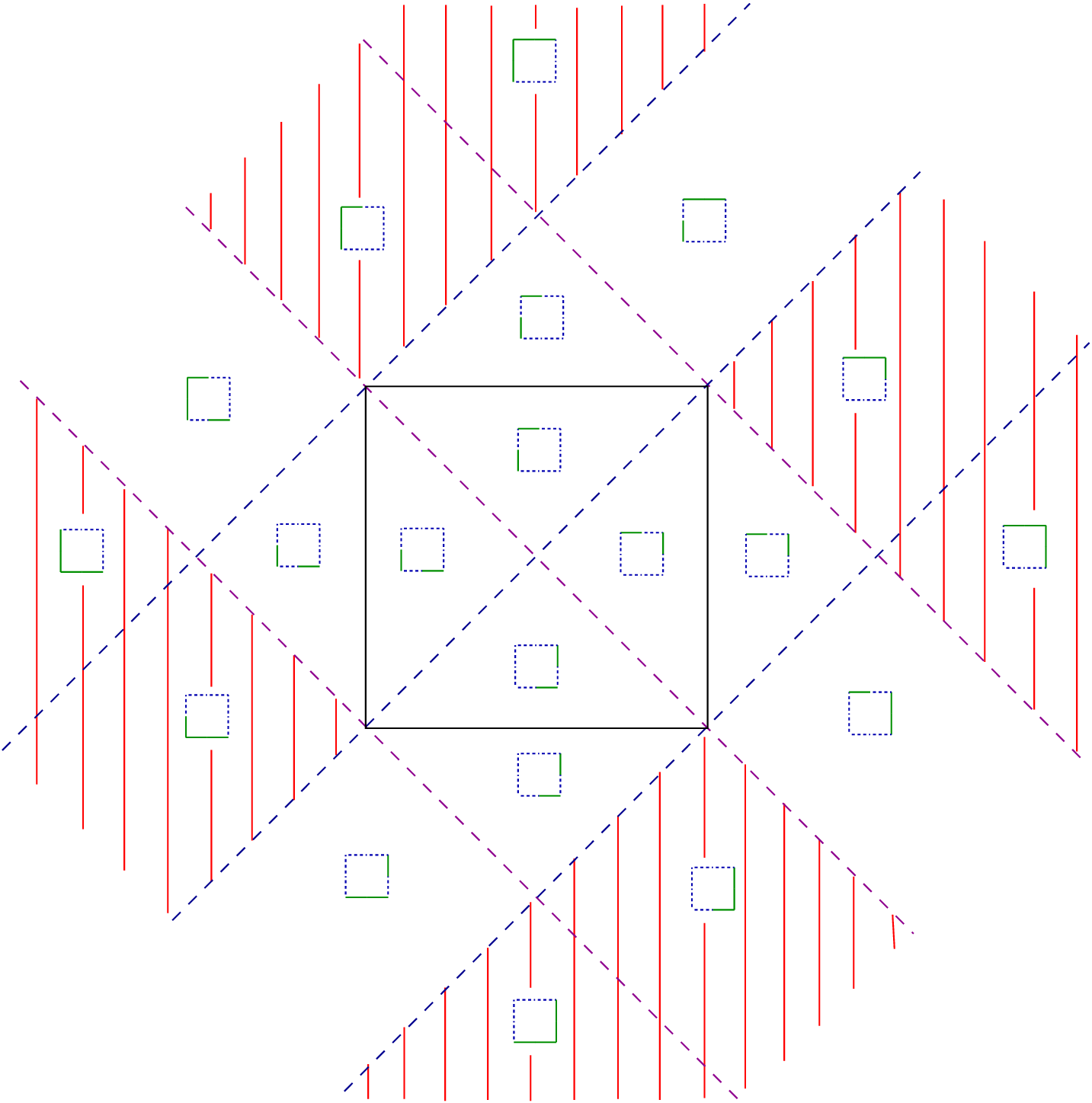,scale=0.18}
\qquad
\epsfig{figure=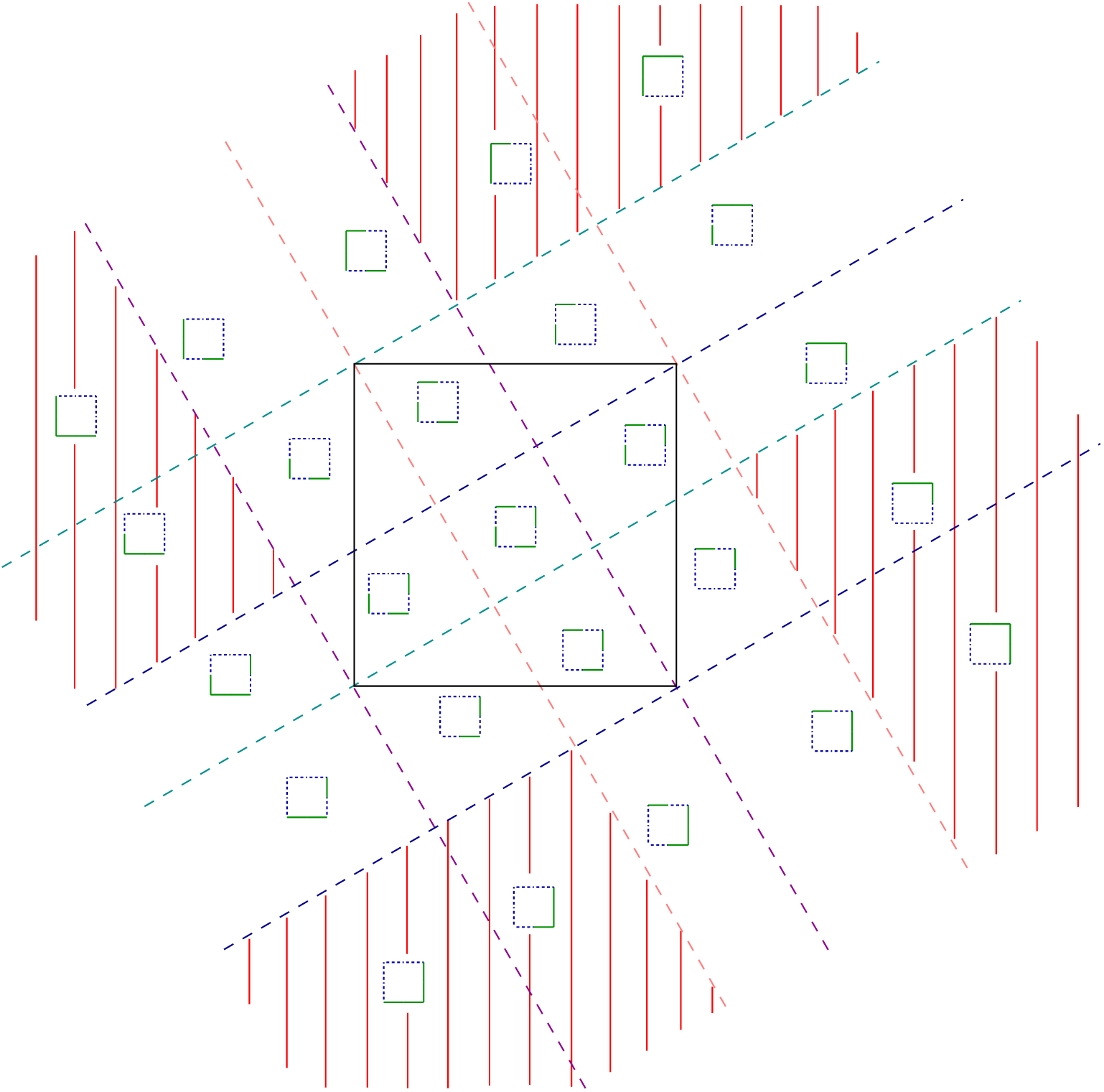,scale=0.18}
}
\caption{Shape of boundary data with respect to $\x_0$ (from left to right,
cases 1,2,3).}
\end{figure}

\subsection{Numerical results}
We perform numerical experiments by considering the following case
$$\Omega =(0,1)^2\, , \quad \partial \Omega_D =
\bigl( \{ 0 \} \times \bigl[ 0,\frac12 \bigr] \bigr) \cup
\bigl( \bigl[ 0,1\bigr] \times \{ 0 \} \bigr) \, ,\quad
\partial \Omega_N =\partial \Omega \setminus \partial \Omega_D \, ,$$
and using above vector field
$$ \m (\x ) =R_\theta (\x -\x_0 ) \, .$$
We only consider the case of a linear feedback. Let us write down the problem.
\begin{equation*}
\left\{ 
\begin{array}{l}
u^{\prime \prime }-\Delta u=0 \\ 
u=0 \\ 
\partial _{\nu }u=-\m .\bn \, u^{\prime } \\ 
u(0)=u_{0} \\ 
u^{^{\prime }}(0)=u_{1} 
\end{array}
\right. \left. 
\begin{array}{l}
\text{in }\Omega \times \mathbb{R}_{+}^{\ast } \, , \\ 
\text{on }\partial \Omega _{D}\times \mathbb{R}_{+}^{\ast } \, , \\ 
\text{on }\partial \Omega _{N}\times \mathbb{R}_{+}^{\ast }\, , \\ 
\text{in }\Omega \, , \\ 
\text{in }\Omega \, .
\end{array}
\right.
\end{equation*}
We will investigate cases when $\theta $ varies in  $[0,\arctan(2)]$.
A particular case is given in Figure 8.
\begin{figure}[htb]
\centerline{
\begin{picture}(0,0)
\includegraphics{carre.pstex}
\end{picture}
\setlength{\unitlength}{4144sp}
\begingroup\makeatletter\ifx\SetFigFont\undefined
\gdef\SetFigFont#1#2#3#4#5{
  \reset@font\fontsize{#1}{#2pt}
  \fontfamily{#3}\fontseries{#4}\fontshape{#5}
  \selectfont}
\fi\endgroup
\begin{picture}(2170,1761)(64,-1660)
\put(1450,-699){\makebox(0,0)[lb]
{\smash{{\color[rgb]{0,0,0}$\Omega $}}}}
\put(2200,-475){\makebox(0,0)[lb]
{\smash{{\color[rgb]{0,0,0}$\partial \Omega_N $}}}}
\put(353,-622){\makebox(0,0)[lb]
{\smash{{\color[rgb]{0,0,0}$\theta $}}}}
\put(79,-1554){\makebox(0,0)[lb]
{\smash{{\color[rgb]{0,0,0}${\mathcal D}_\theta $}}}}
\put(1284,-1500){\makebox(0,0)[lb]
{\smash{{\color[rgb]{0,0,0}$\partial \Omega_D $}}}}
\put(2000,-1460){\makebox(0,0)[lb]{\smash{$\mathbf b$
}}}
\put(2000, 0){\makebox(0,0)[lb]{\smash{$\mathbf c$
}}}
\put(750,-1460){\makebox(0,0)[lb]{\smash{$\mathbf a$
}}}
\put(750, 0){\makebox(0,0)[lb]{\smash{$\mathbf d$
}}}
\put(900,-700){\makebox(0,0)[lb]{\smash{${\boldsymbol \alpha }$
}}}
\end{picture}
}
\caption{When $\x_0 $ belongs to ${\mathcal D}_\theta $, geometrical condition
\eqref{S2} is satisfied at ${\boldsymbol \alpha }$.}
\end{figure}

\noindent Our aim here is to study numerically the variations of the speed of 
stabilization with respect to the position of $\x_0$ and the value of $\theta $.\\
To this end, we have built a finite differences scheme (in space). This leads to
a linear second order differential equation
\begin{equation}\label{EDO}
U''+BU'+KU=0\, ,
\end{equation}
where $B$ is the feedback matrix and $-K$ is the discretized Laplace
operator.\\
Let us define $\,\displaystyle{V=K^{1/2}U} $. Above differential equation can 
be rewritten as follows
\[ \left( \begin{array}{l}
     V\\
     U'
   \end{array} \right)' = \left( \begin{array}{c}
     0\\
     - K^{1/2}
   \end{array} \begin{array}{l}
     K^{1/2}\\
     - B
   \end{array} \right) \left( \begin{array}{l}
     V\\
     U'
   \end{array} \right) \]
and the energy function can be approximated by
$\, \displaystyle{\frac{1}{2}(\langle U,KU\rangle+\|U'\|^2)=
\frac{1}{2}(\|V\|^2+\|U'\|^2)}$.\\
The decreasing rate is given by the highest eigenvalue of above matrix.
Results of our computations are shown in Figure 9 where we built the decreasing
rate as a function depending on $\theta $ and the position of $\x_0 $ represented
by the abscissa $\lambda $ along ${\mathcal D}_\theta $.

\begin{figure}[htb]
\centerline{\epsfig{figure=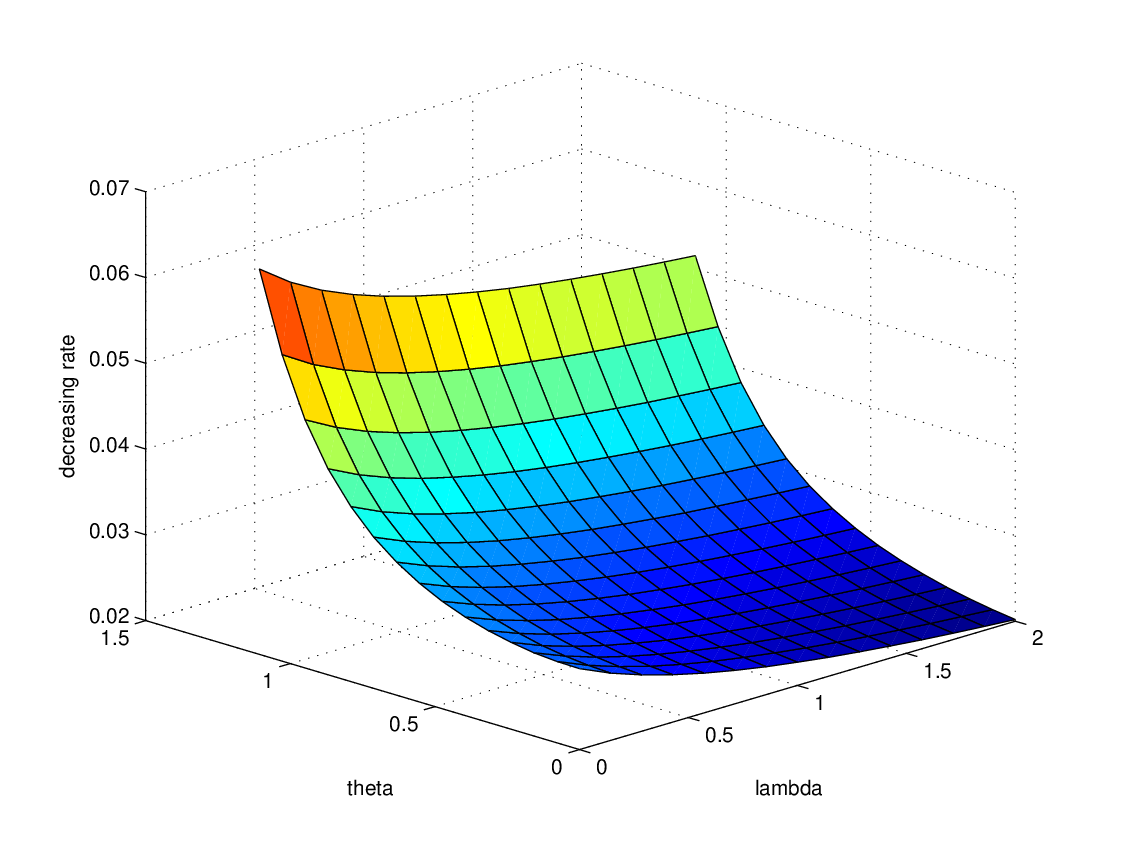,scale=0.7}}
\caption{Dependance of the decreasing rate with respect to 
$\theta$, $\lambda$.}
\end{figure}

\noindent It can be observed that in this case, the decreasing rate is
increasing with $\theta $ and the best position for $\x_0 $ is the
origin of half-line ${\mathcal D}_\theta $. 

\medskip

\begin{small}
\noindent {\bf Acknowledgments}. This work is partly supported by 
FONDECYT 1061263 and ECOS-CONICYT CO4E08.
\end{small}

\end{document}